\documentclass[oneside,a4paper]{article}
\usepackage[utf8]{inputenc}
\usepackage[width=16cm,height=21cm]{geometry}
\usepackage{enumitem}
\usepackage[many]{tcolorbox}
\usepackage{graphicx}
\usepackage{pdfpages}
\usepackage{amssymb,amsmath,amsthm,dsfont}
\usepackage{authblk}
\usepackage{nicefrac}
\usepackage{stmaryrd}
\usepackage{empheq}
\usepackage{amsmath}
\usepackage{todonotes}
\usepackage{url}
\usepackage{hyperref}
\usepackage[capitalise]{cleveref}
\usepackage[mathlines]{lineno}

\newcommand*\linenomathpatch[1]{%
  \cspreto{#1}{\linenomath}%
  \cspreto{#1*}{\linenomath}%
  \csappto{end#1}{\endlinenomath}%
  \csappto{end#1*}{\endlinenomath}%
}
\newcommand*\linenomathpatchAMS[1]{%
  \cspreto{#1}{\linenomathAMS}%
  \cspreto{#1*}{\linenomathAMS}%
  \csappto{end#1}{\endlinenomath}%
  \csappto{end#1*}{\endlinenomath}%
}

\expandafter\ifx\linenomath\linenomathWithnumbers
  \let\linenomathAMS\linenomathWithnumbers
  \patchcmd\linenomathAMS{\advance\postdisplaypenalty\linenopenalty}{}{}{}
\else
  \let\linenomathAMS\linenomathNonumbers
\fi

\linenomathpatch{equation}
\linenomathpatchAMS{gather}
\linenomathpatchAMS{multline}
\linenomathpatchAMS{align}
\linenomathpatchAMS{alignat}
\linenomathpatchAMS{flalign}

\makeatletter
\patchcmd{\mmeasure@}{\measuring@true}{
  \measuring@true
  \ifnum-\linenopenaltypar>\interdisplaylinepenalty
    \advance\interdisplaylinepenalty-\linenopenalty
  \fi
  }{}{}
\makeatother

\newcommand\nnfootnote[1]{%
  \begin{NoHyper}
  \renewcommand\thefootnote{}\footnote{#1}%
  \addtocounter{footnote}{-1}%
  \end{NoHyper}
}

\newcommand{\R}{\mathbb{R}}
\newcommand{\N}{\mathbb{N}}
\newcommand{\Z}{\mathbb{Z}}

\newcommand{\e}{\varepsilon}
\newcommand{\di}[1]{\,\mathrm{d}#1}
\newcommand{\dive}{\operatorname{div}}
\newcommand{\ddt}{\frac{\operatorname{d}}{\operatorname{d}t}}

\newcommand{\ID}{\operatorname{Id}}
\newcommand{\trace}{\operatorname{tr}}
\newcommand{\interior}{\operatorname{int}}
\newcommand{\twosc}{\stackrel{2}{\rightharpoonup}}
\newcommand{\stwosc}{\stackrel{2}{\rightarrow}}

\newcommand{\p}{\varphi}

\newtheorem{theorem}{Theorem}[section]
\newtheorem{lemma}[theorem]{Lemma}
\newtheorem{corollary}{Corollary}[theorem]

\theoremstyle{definition}
\newtheorem{definition}[theorem]{Definition}

\theoremstyle{remark}
\newtheorem{remark}[theorem]{Remark}

\usepackage{color}

\newtcolorbox{mybox}[1]{%
    tikznode boxed title,
    enhanced,
    arc=0mm,
    interior style={white},
    attach boxed title to top left= {yshift=-\tcboxedtitleheight/2-0.05cm, xshift=0.7cm},
    fonttitle=\small\bfseries,
    colbacktitle=white,coltitle=black,
    boxed title style={size=small,colframe=white,boxrule=0pt},
    title={#1}}

\begin{document}
\title{Thermo-Elasticity Problems with Evolving Microstructures}



\date{\today}
\author[$\star,1$]{Michael Eden}
\author[$\dagger$]{Adrian Muntean}

\affil[$\star$]{Faculty of Mathematics, University Regensburg, Germany}
\affil[$\dagger$]{Department of Mathematics and Computer Science, Karlstad University, Sweden}
\maketitle

\begin{abstract}
We consider the mathematical analysis and homogenization of a moving boundary problem posed for a highly heterogeneous, periodically perforated domain.
More specifically, we are looking at a one-phase thermo-elasticity system with phase transformations where small inclusions, initially periodically distributed, are growing or shrinking based on a kinetic under-cooling-type law and where surface stresses are created based on the curvature of the phase interface.
This growth is assumed to be uniform in each individual cell of the perforated domain.
After transforming to the initial reference configuration (utilizing the Hanzawa transformation), we use the contraction mapping principle to show the existence of a unique solution for a possibly small but $\e$ independent time interval ($\e$ is here the scale of heterogeneity).

In the homogenization limit, we recover a macroscopic thermo-elasticity problem which is strongly non-linearly coupled (via an internal parameter called height function) to local changes in geometry.
As a direct by-product of the mathematical analysis work, we present an alternative equivalent formulation which lends itself to an effective pre-computing strategy that is very much needed as the limit problem is computationally expensive.
\end{abstract}

\maketitle
\setcounter{footnote}{1} 
\nnfootnote{$^1$Corresponding author, Email: \href{mailto:michael.eden@ur.de}{michael.eden@ur.de}}

\section{Introduction}\label{sec:introduction}
Phase transformations are complex processes in which different phases grow or shrink at the expense of other phases.
Notable examples are liquid-solid transitions as in the water-ice system or solid-solid transformations as they are happening in alloy steels with, e.g., the passage from austenite to pearlite phases.
Especially in the cases of solidification or solid-solid transformations, these dynamics often occur on the microscale and can lead to complex patterns and microstructures that influence the macroscopic properties of the system \cite{bhadeshia1999,EKK02,pereloma2012}.
Although these transformations are usually driven by temperature, mechanical effects often play a major role in the overall process as potentially different phase densities and localized surface stresses are intrinsically involved \cite{denis1985,mahnken2012}.
As critical geometric changes happen on a length scale much smaller than the characteristic scales observable in laboratory experiments, it is crucially important to derive effective models which are able to explain how measurable material parameters change depending on the microstructure evolution.

In this work, we investigate a prototypica eat conduction in an active medium is coupled one way to a quasi-stationary elasticity system.
The underlying geometry is a connected, perforated domain that is initially periodic.
The small periodicity parameter $\e\ll1$ indicates where the perforations (holes) are allowed to grow or shrink according to the local temperature average via kinetic undercooling (see, e.g., \cite{Back14,PrussSimonettZacher13}) but without surface tension.
The curvature of the phase interface and the resulting surface stresses are accounted for in the mechanics part of the problem.
We use the Hanzawa transformation framework \cite{Hanzawa1981,PrussSimonettZacher13}, where phase growth is characterized by a height function that tracks the changing geometry.
This framework, which sometimes is also called \textit{direct mapping method} (e.g.,\cite{Prokert1999}), allows for a quite general treatment of moving boundary problems.
After transforming the problem into fixed coordinates, we are able to show that there exists a unique weak solution to the moving boundary problem for a small but $\e$-independent time interval.
We then conduct a rigorous limit process $\e\to0$ via two-scale convergence to arrive at a homogenized two-scale model, where the changes at the microscopic level are accounted for via the coefficients of the limit partial differential equation. 

In most cases, it is difficult to establish rigorous homogenization results for moving boundary problems.
The main reason is the necessity for results of higher regularity with delicate control in terms of the small parameter $\e$.
For specific problems, however, a couple of results are available; we point out obstacle and flame propagation problems; see, for instance, \cite{karakhanyan2016,karakhanyan2014}, and, respectively, \cite{caffarelli2006}.
Related results on the homogenization of Stefan-type problems are given in \cite{rodrigues1982}, where a one-phase Stefan problem with rapidly oscillating coefficients was considered, and in \cite{visintin2007homogenization}, where a doubly nonlinear two-phase problem is considered.
For moving boundary problems in which inclusions of a perforated domain are allowed to grow or shrink, most results are restricted either to the linear case, where the evolution is prescribed (e.g., \cite{EKK02,eden_homogenization_2019, Gahn21}), or to phase field approximations, where the geometric changes are not explicitly tracked (e.g., \cite{hoepker16_diss,MSA15,MBW08}).
Steps towards handling fully coupled moving boundary problems posed in perforated domains have been done recently. We mention in particular \cite{Wiedemann}, where the authors investigated the homogenization of a reaction-diffusion-precipitation system, where the evolution of microstructures is explicitly parameterized via a radius function describing the growth/shrinkage of balls.
A similar model was investigated and homogenized in \cite{GahnPop23}.
Both these works feature a somewhat restricted growth rate with hard-coded growth stops, i.e., there is a maximal (minimal) radius from where no further growth (shrinkage) is possible.

Compared to the existing literature, our work brings several novelties:

\begin{itemize}
    \item \textbf{No stopping criteria.} There are no stopping criteria encoded in the growth functions.
    As a consequence, the existence result is only local in time.
    In particular, showing the existence of a $\e$-independent time interval (which makes it possible to consider $\e\to0$) is more involved.
    \item \textbf{More general geometry.} Although we still assume a uniform growth at the level of the cells, we are considering a more general geometric setup compared to the growing or shrinking balls.
    In our setup, any $C^3$-regular inclusions can be considered.
    This is done by characterizing the evolution by a height function and employing the \textit{Hanzawa transformation} \cite{Pruss2016}.
    Since this transformation can also be used to describe a non-uniform cell evolution, this is also an important step to this more general case.
    \item \textbf{Interfacial stresses.} Related works consider the case of precipitation of chemical substances on the surfaces inside porous media (cf.~\cite{Gahn21,GahnPop23,Muntean2020}).
    Although these models are similar to phase transition models, there are also important differences.
    For one, inherently geometric quantities like the curvature (to account for surface stresses) have to be tracked within the problem.
\end{itemize}

This article is structured as follows. In \cref{sec:setting}, we introduce the setup of the moving geometry and the mathematical model equations of our system.
This is followed by \cref{sec:weak_form}, where we list and discuss our mathematical assumptions. We also present here our concept of weak solution.
In \cref{sec:interface_tracking}, we investigate the mathematical analysis of the $\e$-dependent problem and show that there is a unique local-in-time weak solution (see \cref{thm:main_existence}).
Next, \cref{sec:homogenization} is dedicated to the derivation of the homogenization limit $\e\to0$ (\cref{thm:homogenization_limit}).
Finally, in \cref{sec:discussion}, we discuss our results and point to future work relevant for this setting.
\section{Setting and Preliminaries}\label{sec:setting}
In this section, we provide the specific geometric setup, present the mathematical model, and collect assumptions on the coefficients and data.
In \cref{ssec:Hanzawa}, we introduce the Hanzawa transform, which we use to track geometric changes, and we state some important auxiliary results. 
We also present some important lemmas that we rely on in the analysis of our problem in \cref{ssec:auxiliary_results}.
In general, we use $C>0$ to denote any generic constant whose precise value might change even from line to line, but is always independent of $\e$ and $T$, as well as, sometimes, other parameters that are explicitly pointed out.

\begin{figure}[ht]
\centering
\begin{tikzpicture}

  \node[anchor=south west, inner sep=0] (img1) at (0,0)
    {\includegraphics[width=0.3\textwidth]{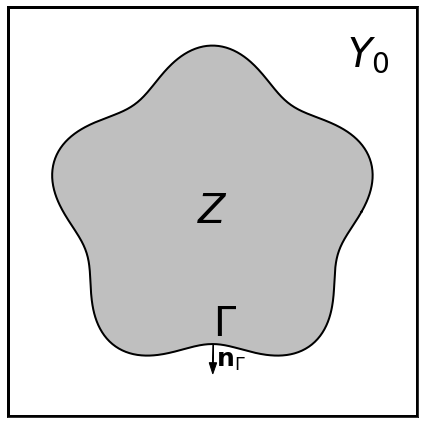}};

  \node[anchor=south west, inner sep=0] (img2) at (7,0)
    {\includegraphics[width=0.3\textwidth]{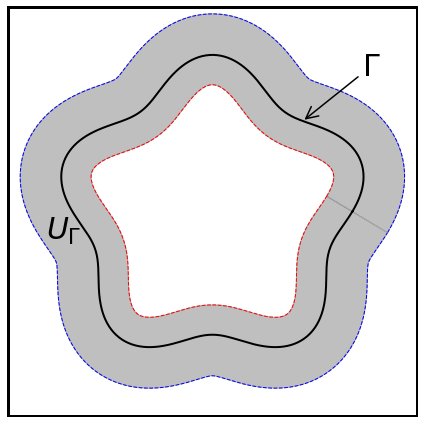}};



\end{tikzpicture}
  \caption{Left: The geometry of the unit cell with the inclusion $Z$.
  Right: Visualization of a tubular neighborhood $U_\Gamma$ (grey area) introduced in \cref{ssec:Hanzawa} inside which the interface is allowed to grow.
  Please note that these are cross sections as we are working in 3D.}
  \label{figure:unit_cell}
\end{figure}

First, let $S=(0,T)$, $T>0$, represent the time interval of interest and let a bounded Lipschitz domain $\Omega\subset\R^3$ represent the spatial domain.
For technical reasons, we assume that $\Omega$ is a finite union of axis-parallel cubes with corner coordinates in $\mathbb{Z}^{3}$.
We denote the outer normal vector of $\Omega$ by $\nu=\nu(x)$.
Now, let $Y=(0,1)^3$ denote the standard unit cell, and let $Z\subset Y$ be a $C^3$ regular domain with $\overline{Z}\subset Y$. 
Removing $\overline{Z}$ from $Y$, we get the Lipschitz domain $Y_0=Y\setminus\overline{Z}$, whose external and internal boundaries are given by $\partial Y$ and $\Gamma:=\partial Z$, respectively.
With $n_\Gamma=n_\Gamma(\gamma)$, $\gamma\in\Gamma$, we denote the normal vector of $\Gamma$ pointing inside of $Y_0$.
See \cref{figure:unit_cell} for a representation of this set-up.

Now, let $\e_0>0$ be chosen such that $\Omega$ can be perfectly tiled with $\e_0Y$ cells and set $\e_n=2^{-n}\e_0$.
For readability, we will omit the $n$ in $\e_n$ and interpret $\e\to0$ as $\e_n\to0$ for $n\to\infty$.
We introduce the finite index set $\mathcal{I}_\e=\{k\in\mathbb{Z}^3\ : \ \e(Y_0+k)\subset\Omega_\e\}$ and use it to index individual cells $Y_{\e,k}:=\e(Y_0+k)$ and their internal interfaces through $\Gamma_{\e,k}:=\e(\Gamma+k)$ for all $k\in \mathcal{I}_\e$.
Please note that $|Y_{\e,k}|=\e^{3}|Y_0|$ and $|\Gamma_{\e,k}|=\e^2|\Gamma|$.
From this we built up the initial perforated domain and its internal interface
    \[
    \Omega_\e:=\interior\left(\bigcup_{k\in\mathcal{I}_\e}\e(\overline{Y_0}+k)\right),
    \quad \Gamma_\e:=\Omega\cap \left(\bigcup_{k\in\mathcal{I}_\e}\e(\Gamma+k)\right).
    \]

\begin{figure}[ht]
\centering
\begin{tikzpicture}

  \node[anchor=south west, inner sep=0] (img1) at (0,0)
    {\includegraphics[width=0.4\textwidth]{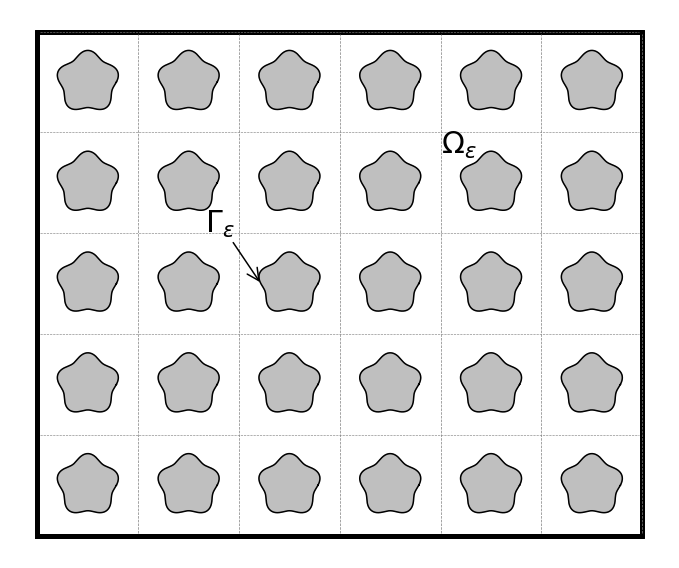}};

  \node[anchor=south west, inner sep=0] (img2) at (7,0)
    {\includegraphics[width=0.4\textwidth]{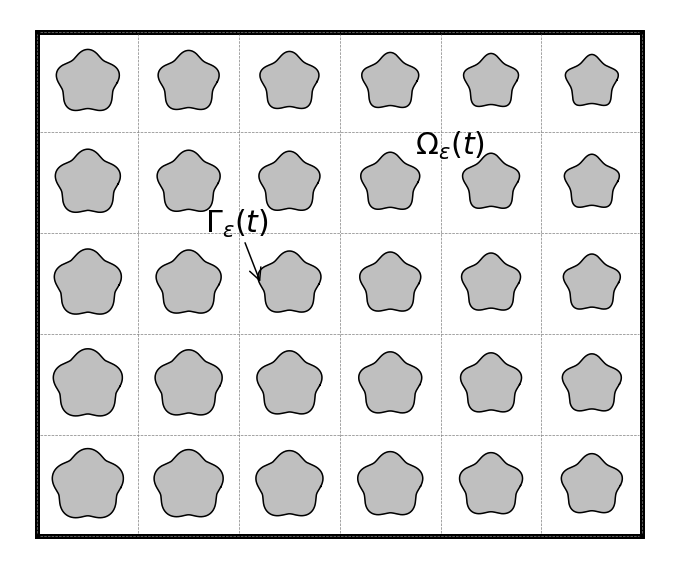}};
\end{tikzpicture}
\caption{Left: A two dimensional cross section of the initial periodic geometry. Right: The changed geometry of this cross section at some time $t>0$ with growth and shrinkage in different parts of the domain.
Please note that these are cross sections as we are working in 3D.}
\label{figure:domain_evolution}
\end{figure}

We introduce time-dependent domains $\Omega_\e(t)$ and $\Gamma_\e(t)$ for $t\in(0,T)$, see \cref{figure:domain_evolution}.
The normal vector at $\Gamma_\e(t)$ that points inside $\Omega_\e(t)$, is denoted by $n_\e(t,\gamma)$ for all $\gamma\in\Gamma_\e(t)$.
We assume that the outer boundary $\partial\Omega$ does not change in time and introduce the space-time sets
\[
Q_\e=\bigcup_{t\in S}\{t\}\times\Omega_\e(t),\qquad \Sigma_\e=\bigcup_{t\in S}\{t\}\times\Gamma_\e(t).
\]

For $t\in S$, and $x\in\Omega_\e(t)$, let $u_\e=u_\e(t,x)$ denote the deformation and $\theta_\e=\theta_\e(t,x)$ the temperature.
The bulk equations for the one-phase thermoelasticity problem are given as (for more details on the modeling, we refer to~\cite{B56,K79})
\begin{subequations}\label{p:full_problem_moving}
\begin{alignat}{2}
	-\dive(\mathcal{C} e(u_\e)-\alpha\theta_\e\mathrm{I}_3)&=f_\e&\quad&\text{in}\ \ Q_\e,\label{p:full_problem_moving:1}\\
	\rho c\partial_t\theta_\e-\dive(K\nabla\theta_\e)&=g_\e&\quad&\text{in}\ \ Q_\e,\label{p:full_problem_moving:3}
\end{alignat}
Here, $\mathcal{C}\in\R^{3\times3\times3\times3}$ is the \emph{stiffness} tensor of elasticity, $\alpha>0$ the \emph{thermal expansion} coefficient, $\rho>0$ the \emph{mass density}, $c>0$ the \emph{heat capacity}, $K\in\R^{3\times3}$ the \emph{thermal conductivity}, and $f_\e$, $g_\e$ are volume source densities for stresses and heat, respectively.
In addition, $e(v)=1/2(\nabla v+\nabla v^T)$ denotes the linearized strain tensor and $\mathrm{I}_3$ the identity matrix.
Here, we are neglecting the heat dissipation term, $\gamma\dive \partial_tu_\e$ with \emph{dissipation coefficient} $\gamma>0$, as is often done in similar problems \cite{Wan99}.

The phase transformation and as a consequence the geometric changes are assumed to be driven by the local average of the temperature via the law of kinetic undercooling:
\begin{equation}\label{p:full_problem_moving:4}
v_{\e,k}=\frac{1}{|\Gamma_{\e,k}(t)|}\int_{\Gamma_{\e,k}(t)}\theta_\e(t,x)\di{\sigma}.
\end{equation}
Here, $v_{\e,k}\colon S\to\R$ for $k\in \mathcal{I}_\e$ denotes the normal velocity (in direction of $n_\Gamma$) of the interfaces.
Therefore, positive average temperatures at the interface $\Gamma_{\e,k}$ lead to cell shrinkage $\e(Y_0+k)$ and, correspondingly, negative averages to growth.

The phase transformation might induce additional surface forces and we assume these forces to be proportional to the mean curvature of the interface:
\begin{equation}\label{p:full_problem_moving:6}
	-(\mathcal{C}e(u_\e)-\alpha\theta_\e\mathrm{I}_3)n_\e=-\e^2\kappa_\e\sigma_0n_\e\quad\text{on}\ \ \Sigma_\e.
\end{equation}
Here, $\sigma_0>0$ is the coefficient of surface tension and $\kappa_\e$ is the mean curvature of the interface.
The scaling via $\e^2$ counters the effects of both the interface surface area and the curvature itself, note that $|\Gamma_\e|,|\kappa_{\e}|\sim\e^{-1}$.

In a similar way, the release of heat across the interface is assumed to be given via the \emph{latent heat} $L\in\R$:
\begin{equation}\label{p:full_problem_moving:7}
	-K\nabla\theta_\e\cdot n_\e=\e Lv_\e\quad\text{in}\ \ \Sigma_\e,
\end{equation}
where $v_\e$ denotes the normal velocity of the interface.
Here, the scaling via $\e$ counters the effect of the interface surface area.
More complex interface conditions than equations~\eqref{p:full_problem_moving:6}, \eqref{p:full_problem_moving:6} arise if the interface is allowed to be thermodynamically active thereby requiring us to formulate separate balance equations for surface stress and surface heat, we refer to~\cite{WB15}.

Finally, we pose homogeneous \emph{Dirichlet conditions} for the momentum equation and homogeneous \emph{Neumann conditions} for the heat equation as well as initial conditions for the temperature:
\begin{alignat}{2}
	u_\e&=0&\quad&\text{on}\ \ S\times\partial\Omega,\\
	-K\nabla\theta_\e\cdot\nu&=0&\quad&\text{on}\ \ S\times\partial\Omega,\\
	\theta_\e(0)&=\theta_{0\e}&\quad&\text{on}\ \ \Omega_\e,\label{p:full_problem_moving:9}
\end{alignat}
where $\theta_{0\e}$ is some (possibly highly heterogeneous) initial temperature distribution.
\end{subequations}

\begin{remark}[The averaging assumption]\label{rem:averaging_assumption}
    Please note that the averaging condition \cref{p:full_problem_moving:4} is chosen for mathematical reasons to avoid local velocity gradients.
    In principle, local gradients can be captured with the Hanzawa transformation, but it is not yet clear how to incorporate those into homogenization settings.
    We note that this approach with local averaging is a relatively common approach to tackle such problems, cf.~\cite{Gahn21,Wiedemann}.
    The main obstacle in treating the case without averaging, i.e, replacing \cref{p:full_problem_moving:4} with $v_{\e,k}=\theta_{\e,k}$, are the additional $\e$-uniform estimates needed in the analysis.
    With the averaging assumption, standard energy estimates together with $L^\infty$-bounds are sufficient (cf.~\cref{lem:reference_apriori}).
    Without it, additional estimates in the form of 
    \[
    \|\nabla\theta_\e\|_{L^\infty(\Omega_\e(t))}+\e\|\nabla^2\theta_\e\|_{L^\infty(\Omega_\e(t))}\leq C
    \]
    with $C$ independent of $\e$ are necessary; see for example \cite{Chen92,eden_homogenization_2019}.
    
    We note that due to the connectedness of $\Omega_\e(t)$, we expect $\theta_\e$ to be strongly converging and, as a result, to be independent of the micro variable in the homogenization limit.
    For that reason, we expect
    \[
    \left(\theta_\e(t,\gamma)-\frac{1}{|\Gamma_{\e,k}(t)|}\int_{\Gamma_{\e,k}(t)}\theta_\e(t,\gamma)\di{\sigma}\right)\to0\quad(\e\to0)
    \]
    almost everywhere on $\Gamma_{\e}(t)$.
    This indicates that this modeling assumption is also physically reasonable in this specific case.
\end{remark}

A comment regarding the notation:
We do not distinguish between the corresponding norms for scalar, vector, or matrix valued spaces, e.g., we write $\|\cdot\|_{L^q(\Omega)}$ for functions in $L^q(\Omega)$, $L^q(\Omega)^3$, as well as $L^q(\Omega)^{3\times3}$.
In each case, the norm is understood as the Bochner norm, with the Euclidean norm used for vector-valued functions and either the Frobenius or operator norm for matrices, depending on the context.
For example, for $\Psi\in L^2(\Omega)^{3\times3}$, we have
\[
\|\Psi\|_{L^2(\Omega)}:=
\left(\int_{\Omega}|\Psi(x)|^2\di{x}\right)^\frac{1}{2}
\]
where $|\cdot|^2$ denotes the norm on $\R^{3\times3}$ (either Frobenius or operator norm).
Since $|A|_{Op}\leq |A|_{Fr}\leq \sqrt{3}|A|_{Op}$ for all $A\in\R^{3\times3}$, we can switch between these norms.
Whenever specific estimates (i.e., without generic constants) are claimed (as in \cref{hanzawa_normal_vel} for $F_\e=Ds_\e$ with $\|F_\e\|_{L^\infty(\Omega)}\leq2$), the operator norm is used.

\subsection{Coordinate transform}\label{ssec:Hanzawa}
We start by introducing some rudimentary concepts of differential geometry and point to \cite{Lee2018,Pruss2016} for an in-depth overview and discussion.
As a $C^3$ interface, $\Gamma$ admits a tubular neighborhood $U\subset Y$ as depicted in \cref{figure:unit_cell}.
Moreover, there exist $a_i>0$ ($i=1,2$)\footnote{These are the radii of balls for which the interface satisfies interior and exterior ball conditions.} such that the function
\[
\Lambda\colon\Gamma\times(-a_1,a_2)\to U,\quad \Lambda(\gamma,s)=\gamma+sn_\Gamma(\gamma)
\]
is a $C^2$ diffeomorphism whose inversion is given by
\[
\Lambda^{-1}\colon U\to\Gamma\times(-a_1,a_2),\quad \Lambda^{-1}(y)=(P(y),d(y))
\]
where $P$ denotes the projection operator onto $\Gamma$ and $d$ the signed distance function (positive in $Y_0$ and negative in $Z$).
Next, we introduce the \emph{shape tensor}, or \emph{Weingarten map}, 
\[
L_\Gamma\colon\Gamma\to\R^{3\times3},\quad
L_\Gamma(\gamma)=D\left(\frac{\nabla\rho(\gamma)}{|\nabla\rho(\gamma)|}\right)
\]
where $\rho$ is any $C^3$ level set function for $\Gamma$.\footnote{The shape tensor is independent of the choice of level set, see \cite[Section 2.1]{Pruss2016}.}
Here and in the following, we use $D$ to denote the Jacobian matrix of any vector valued function and note that $Df=(\nabla f)^T$ for any real valued function $f\colon\R^3\to\R$. The Weingarten map is symmetric and its operator norm therefore is given by $|L_\Gamma(\gamma)|=\max_{j=1,2}|\kappa_j|$ where the eigenvalues $\kappa_1,\kappa_2\in\R$ are the principial curvatures of $\Gamma$ at $\gamma\in\Gamma$.\footnote{The third eigenvalue is $0$ as $L_\Gamma n_\Gamma=0$.}
Finally, we set $a^*=\min\{a_1,a_2\}$ and have $0<a^*\leq(2|L_\Gamma|(\gamma)|)^{-1}$ which implies $|L_\Gamma(\gamma)|\leq (2a^*)^{-1}$.
The mean curvature field $\kappa\colon\Gamma\to\R$ of the interface $\Gamma$ is given by
\[
\kappa(\gamma)=\frac{1}{2}\trace L_\Gamma(\gamma)
\]
and we have $\kappa\in C^1(\Gamma)$ since $\Gamma$ is $C^3$.
Finally, we set $a^*=\min\{a_1,a_2\}$.
In the following lemma, we characterize some derivatives that will come up in the analysis.

\begin{lemma}\label{lem:derivatives_geometry}
    The derivatives of $d$ and $P$ are given via ($y\in U$)
    \begin{align*}
    \nabla d(y)&=n_\Gamma(P(y)),\\
    DP(y)&=M(P(y),d(y))\left(\ID-n_\Gamma(P(y))\otimes n_\Gamma(P(y))\right)
    \end{align*}
    where
    \[
    M\colon\Gamma\times(a_1,a_2)\to\R^{3\times3},\quad M(\gamma,s)=(\ID-sL_\Gamma(\gamma))^{-1}.
    \]
    In addition, we have the implicit relation
    \[
    D^2d(y)=-L_\Gamma(P(y))(\ID-d(y)L_\Gamma(P(y)))^{-1}.
    \]
\end{lemma}

\begin{proof}
    The derivatives are given in \cite[Chapter 2, Section 3.1]{Pruss2016}.
    The invertability of $(\ID-sL(\cdot))^{-1}$ over $\Gamma\times(a_1,a_2)$ follows by $2|L(\gamma)|\leq \min\{a_1,a_2\}^{-1}=\nicefrac{1}{a^*}$. 
\end{proof}
Generally speaking, a Hanzawa transform is a coordinate transformation between different interfaces which are normally parametrizable; i.e., where the change can be measured in terms of a height function describing the change of $\Gamma$ inside the tubular neighborhood $U$.
As we are considering uniform growth (or shrinkage) in each individual cell (see \cref{p:full_problem_moving:4}), this transformation is easy to describe.
To that end, let $\chi\colon\R\to[0,1]$ be a smooth cut-off function which satisfies
\begin{alignat*}{2}
(i)&\ \chi(r)=1 \ \text{for}\ r\in(-\frac{a_1}{3},\frac{a_2}{3}), &\qquad (ii)&\ \chi(r)=0 \ \text{for}\ r\notin(-\frac{2a_1}{3},\frac{2a_2}{3}),\\
(iii)&\ \mathrm{sign}(r)\chi'(r)\leq0,&\quad (iv)&\ |\chi'(r)|\leq\frac{4}{a^*}.
\end{alignat*}

\begin{figure}
    \centering
    \includegraphics[width=0.45\textwidth]{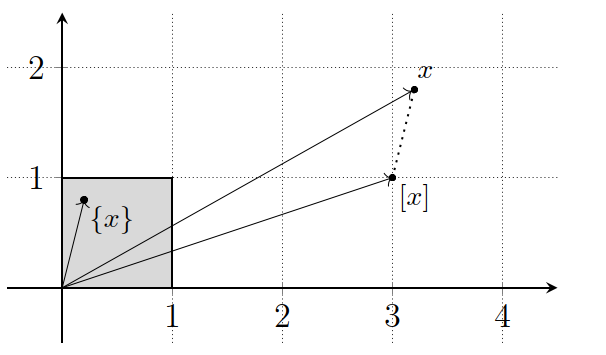}
    \caption{Simple example demonstrating the construction of $[x]$ and $\{x\}$ in $\R^2$}
    \label{fig:unfolding}
\end{figure}

\begin{lemma}\label{lem:hanzawa_stationary}
    Let $h\in(-\nicefrac{a_1}{10},\nicefrac{a_2}{10})$ and $\Gamma(h)=\{\gamma+hn(\gamma)\ : \gamma\in\Gamma\}$.
    Then the $C^2$ transformation $s\colon\overline{Y}\to\overline{Y}$ defined by 
    \begin{equation}\label{hanzawa_transform}
     s(x)=
     x+hn_\Gamma(P(x))\chi\left(d(x)\right)
    \end{equation}
    satisfies $s(\Gamma)=\Gamma(h)$ as well as
    \[
    \|Ds\|_{L^\infty(\Omega)}\leq 2,\quad \|Ds^{-1}\|_{L^\infty(\Omega)}\leq 2.
    \]
    The curvature of $\Sigma$ can be calculated via
    \[
    \kappa\colon(-\frac{a_1}{10},\frac{a_2}{10})\times\Gamma\to\R,\quad \kappa(h,\gamma)=\trace\left[(\ID-hL_\Gamma(\gamma))^{-1}L_\Gamma(\gamma)\right].
    \]
\end{lemma}
\begin{proof}
    Since $\Gamma(h)$ is in the tubular neighborhood of $\Gamma$ and the height parameter $h$ is constant, $s(\Gamma)=\Gamma(h)$ and the $C^2$-regularity of $s_\e$ follows immediately.
    We introduce the function $f\colon[-a_1,a_2]\to[-a_1,a_2]$ with $f(r)=r+h\chi(r)$.\footnote{This function $f$ describes the behavior of $s$ alongside any normal ray.}
    With $f'(r)=1+h\chi'(r)$ and $|\chi'|\leq \nicefrac{4}{a^*}$ injectivity of $f$ follows for all $h\in(-\nicefrac{a_1}{10},\nicefrac{a_2}{10})$.
    In addition, as $f$ is continuous and $f(-a_1)=-a_1$ and $f(a_2)=a_2$, $f$ is also onto.
    Now, assume there are two points $x_1,x_2\in U$ such that $s(x_1)=s(x_2)$ which implies $P(x_1)=P(x_2)$ (points are only moved in normal direction).
    With $d(s(x_1))=d(s(x_2))$, we find that
    \[
    d(x_1)+h\chi(d(x_1))=d(x_2)+h\chi(d(x_2)).
    \]
    Since $f$ is injective, we have $x_1=x_2$; hence injectivity of $s$. From here, the regularity of $s^{-1}$ and the bounds for the Jacobians can be argued via direct calculations.
    We refer to \cite[Lemma 2.9]{eden2018homogenization} where this is done in detail.
    The curvature identity is a special case of \cite[Section  1.2]{Pruss2016}.
\end{proof}
Note that the range of possible $h$ can be extended (too, e.g., $(-\nicefrac{a_1}{2},\nicefrac{a_2}{2})$) at the cost of weaker estimates for $\|Ds^{-1}\|_{L^\infty(\Omega)}$. 
As a next step, we extend this transformation to the rescaled, $\e$-periodic initial geometry.
To that end, we introduce $\e$-periodic interface and its corresponding tubular neighborhood
\[
\Gamma_{\e}:=\Omega\cap \left(\bigcup_{k\in\mathbb{Z}^3}\e(\Gamma+k)\right),\quad
U_{\e}:=\Omega\cap \left(\bigcup_{k\in\mathbb{Z}^3}\e(U+k)\right).
\]
By a simple scaling argument, this leads to the $C^2$-diffeomorphism
\[
\Lambda_{\e}\colon\Gamma_{\e}\times (-\e a_1,\e a_2)\to U_{\e}
\]
with the corresponding inverse
\[
\Lambda_{\e}^{-1}\colon U_{\e}\to\Gamma_{\e}\times(-\e a_1,\e a_2),\quad \Lambda_{\e}^{-1}(x)=(P_{\e}(x),d_{\e}(x))
\]
where $P_{\e}$ denotes the projection operator onto $\Gamma_{\e}$ and $d_{\e}$ the signed distance function.
Please note that, by construction, $P_{\e}$ and $\Gamma_{\e}$ are directly related to the $Y$-periodic functions $P$ and $\Gamma$ via
\[
P_{\e}(x)=\e\left[\frac{x}{\e}\right]+P\left(\e\left\{\frac{x}{\e}\right\}\right),\quad d_{\e}(x)=\e d\left(\left\{\frac{x}{\e}\right\}\right).
\]
Here, $\left[\cdot\right]\colon\R^n\to\Z^n$ and $\left\{\cdot\right\}\colon\R^n\to Y$ are defined via (see also \cref{fig:unfolding})
\[
\left[x\right]=k\ \text{such that}\ x-[x]\in[0,1)^n,\quad \left\{x\right\}=x-[x].
\]

Due to this periodicity, the results from \cref{lem:hanzawa_stationary} can be directly transferred:

\begin{corollary}\label{cor:hanzawa_stationary_eps}
    Let $h\colon\mathcal{I}_\e\to(-\nicefrac{a_1}{10},\nicefrac{a_2}{10})$ and let
    \[
    \Gamma_\e(h)=\bigcup_{k\in\mathcal{I}_\e}\{\gamma+\e h(k)n_{\e}(\gamma)\ : \gamma\in\e(\Gamma+k)\}.
    \]
    Then the $C^2$ transformation $s_\e\colon\overline{\Omega}\to\overline{\Omega}$ defined by 
    \begin{equation}\label{hanzawa_transform_eps}
     s_\e(x)=
     x+\e h(P_{\e}(x))\,n_{\e}(P_{\e}(x))\chi\left(\frac{d_{\e}(x)}{\e}\right)
    \end{equation}
    satisfies $s(\Gamma_{\e})=\Gamma_\e(h_\e)$ as well as
    \[
    \|Ds_\e\|_{L^\infty(\Omega)}\leq 2,\quad \|Ds_\e^{-1}\|_{L^\infty(\Omega)}\leq 2.
    \]
    For the curvature of $\Gamma_\e(h)$ with respect to the curvature of $\Gamma$, we have
    \[
    \kappa_\e\colon(-\frac{a_1}{2},\frac{a_2}{2})\times\Gamma\to\R,\quad\kappa_\e(h,\gamma)=\e^{-1}\trace\left[(\ID-hL_\Gamma(\gamma))^{-1}L_\Gamma(\gamma)\right].
    \]
\end{corollary}
\begin{proof}
    Taking any $\e$-cell $\e(Y+k)$ for $k\in\mathcal{I}_\e$, we can apply \cref{lem:hanzawa_stationary} to get a transformation $s_{\e,k}\colon \e(\overline{Y}+k)\to\e(\overline{Y}+k)$.
    Now, since $s_{\e,k}=\ID$ in a neighborhood of $\partial Y$ due to the cut-off function $\chi$ vanishing, we can glue these transformations smoothly via $s_{\e}(x):=\sum_{k\in\mathcal{I}_\e}\mathds{1}_{\e(\overline{Y}+k)}(x)\hat{s}_{\e,k}(x)$, where $\hat{s}_{\e,k}$ denotes the zero-extension of $s_{\e,k}$.
    Since the $\e(Y+k)$-cells perfectly tile the domain $\Omega$ for all $\e$, this implies $s_\e\colon\overline{\Omega}\to\overline{\Omega}$ as well as its $C^2$-regularity.

    The curvature identity is just a basic rescaling of the corresponding result in \cref{lem:hanzawa_stationary}.
\end{proof}
To account for the changes in time, let, for $t\in\overline{S}$, $\Gamma_\e(t)\subset\mathcal{U}_{\e}$ and let $h\colon\overline{S}\times\mathcal{I}_\e\to(-\nicefrac{a_1}{10},\nicefrac{a_2}{10})$ be height functions satisfying
\[
\Gamma_\e(t)=\bigcup_{k\in\mathcal{I}_\e}\{\gamma+\e h(t,k)n_{\e}(\gamma)\ : \ \gamma\in\e(\Gamma+k)\}. 
\]
Then, we can apply \cref{cor:hanzawa_stationary_eps} at every point in time to get a parameterized family of transformations $s_\e(t,\cdot)$ whose regularity with respect to time is determined by the time regularity of $h$.
With the normal velocity function $v_\e$, the height functions are explicitly given via $h_\e(t,k)=\int_0^tv_\e(\tau,k)\di{\tau}$.\footnote{Please note that this is only possible since the cell growth is uniform. Without this assumption, the height function is more difficult to characterize, see \cite{Pruss2016}.}

\begin{remark}[Example: Growing and shrinking balls.]
    One relatively simple example of this setup with the Hanzawa transformation is balls which are allowed to grow and shrink.
    To that end, let $\Gamma=\partial B_r(y)$, where $B_r(y)$ denotes the ball of radius $r$ centered at $y=(0.5,0.5)$.
    In this case, we may choose $a_1=a_2=0.25$ and can refer to \cite[Section 2.2.1.]{Wiedemann} for an explicit form of a corresponding transformation.
\end{remark}

\subsection{Additional auxiliary results}\label{ssec:auxiliary_results}
There are some additional auxiliary results that we will rely on in the analysis of the full problem.
To formulate these results, we need to talk about how the geometries change with the height functions.
To that end, for $h\colon\mathcal{I}_\e\to(-\nicefrac{a_1}{10},\nicefrac{a_2}{10})$, we set (as in \Cref{cor:hanzawa_stationary_eps})
\[
\Gamma_\e(h):=\bigcup_{k\in\mathcal{I}_\e}\{\gamma+\e h_{\e}(k)n_{\e}(\gamma)\ : \gamma\in\e(\Gamma+k)\}.
\]
In addition, we take $\Omega_\e(h)\subset\Omega$ to be the corresponding connected domain with internal boundary $\Gamma_\e(h)$, that is, $\partial\Omega_\e(h)=\Gamma_\e(h)\cup\partial\Omega$.

\begin{lemma}[Extension operators]\label{lemma:extension_operator}
Let $h\colon\mathcal{I}_\e\to(-\nicefrac{a_1}{10},\nicefrac{a_2}{10})$.
There is a family of linear extension operators $\mathcal{P}_{\e,h}\colon H^1(\Omega_{\e}(h))^3\to H^1(\Omega)^3$ and a constant $C>0$ (which is independent of $\e$ and $h$) such that
\[
\|\mathcal{P}_{\e,h}\varphi\|_{H^1(\Omega)}\leq C\|\varphi\|_{H^1(\Omega_{\e}(h))},\quad \|e(\mathcal{P}_{\e,h}\varphi)\|_{L^2(\Omega)}\leq C\|e(\varphi)\|_{H^1(\Omega_{\e}(h))}
\]
for all $\varphi\in H^1(\Omega_{\e}(h))^3$.
Here, $e(\phi)=\nicefrac{1}{2}(\nabla \varphi+\nabla\varphi^T)$ denotes the symmetric gradient.
\end{lemma}
\begin{proof}
    For every constant $h$ (i.e., periodic geometry), this is a standard result for perforated domains (we refer to \cite[Section 2.3]{Cioranescu} and \cite[Section 4.3]{hoepker16_diss} for the symmetric gradient); and since we are only looking at inclusions and have no changes in topology or connectivity, we can choose the constant $C$ independent of $h$, i.e., 
    \[
    C:=\sup_{-\nicefrac{a_1}{10}<h<\nicefrac{a_2}{10}}C_h<\infty.
    \]
    That this is indeed possible is shown in full mathematical rigor and in a more general and abstract setup in \cite[Theorem 3.1]{gahn2024}.
    
    Now, for variable $h\colon\mathcal{I}_\e\to(-\nicefrac{a_1}{10},\nicefrac{a_2}{10})$, we can locally extend and glue together in the same way as in the standard periodic case like  \cite[Section 2.3]{Cioranescu}.
\end{proof}

\begin{lemma}[Trace operators]\label{trace_operators}
There is a $C>0$ such that
\begin{itemize}
    \item[$(i)$]\textbf{Trace theorem.}
    \[
    \e\|v\|^2_{L^2(\Gamma_{\e})}\leq C\left[\|v\|^2_{L^2(\Omega_{\e})}+\e^2\|\nabla v\|^2_{L^2(\Omega_{\e})}\right]\quad (v\in H^1(\Omega_{\e})).
    \]
    \item[$(ii)$]\textbf{Trace theorem for an individual cell.}
    \[
    \|v\|_{L^1(\Gamma_{\e,k})}\leq C\left[\frac{1}{\e}\|v\|_{L^1(\e(Y+k))}+\|\nabla v\|_{L^1(\e(Y+k))}\right]\quad (v\in H^1(\e(Y+k))).
    \]
    \item[$(iii)$]\textbf{Trace theorem for bounded functions.}
    \[
     \|v\|_{L^\infty(\Gamma_{\e})}\leq \|v\|_{L^\infty(\Omega_{\e})}\quad (v\in H^1(\Omega_{\e})\cap L^\infty(\Omega_\e)).
    \]
\end{itemize}
\end{lemma}
\begin{proof}
    For $(i)$, we refer to \cite[Lemma 3(b)]{HJ91}.
    For $(ii)$, we rely on a scaling argument.
    We note that $\Gamma_{\e,k}=\e(\Gamma+k)$ and, w.l.o.g., take $k=0$ as $k$ is just a translation.
    We employ the substitution $y\mapsto \tilde{y}=\nicefrac{y}{\e}$ and get
    \[
    \int_{\Gamma_{\e,0}}|v(y)|\di{\sigma}=\e^2\int_{\Gamma}|v(\tilde{y})|\di{\tilde{\sigma}}
    \leq \e^2 C\left(\int_{Y}|v(\tilde{y})|\di{\tilde{y}}+\int_{Y}|\nabla_{\tilde{y}}v(\tilde{y})|\di{\tilde{y}}\right)
    \]
    Substituting back, yields
    \[
    \int_{\Gamma_{\e,0}}|v(y)|\di{\sigma}
    \leq \e^{-1} C\left(\int_{\e Y}|v(y)|\di{y}+\e\int_{\e Y}|\nabla_yv(y)|\di{y}\right).
    \]
    For $(iii)$, since $\tilde{v}:=v+\|v\|_{L^\infty(\Omega_\e)}\in H^1(\Omega_\e)$ is non-negative, it holds $\tilde{v}\geq0$ on $\Gamma_\e$ in the sense of traces.
    Consequently, $v\geq -\|v\|_{L^\infty(\Omega_\e)}$ on $\Gamma_\e$.
    The upper bound $v\leq \|v\|_{L^\infty(\Omega_\e)}$ follows via $v-\|v\|_{L^\infty(\Omega_\e)}$.
\end{proof}
We introduce 
\[
H_{\partial\Omega}(\Omega_\e)=\{u\in H^1(\Omega_\e)\ : \ u=0\ \text{on}\ \partial\Omega\}
\]
where we have the following inequalities:
\begin{lemma}[Inequalities for perforated domains]\label{trace_embedding}
There is a $C>0$ such that
\begin{itemize}
    \item[$(i)$]\textbf{Poincare's inequality.}
    \[
     \|v\|_{L^2(\Omega_{\e})}+\|\nabla v\|_{L^2(\Omega_{\e})}
    	\leq C\|\nabla v\|_{L^2(\Omega_{\e})} \quad (v\in H_{\partial\Omega}(\Omega_\e)^{3}).
    \]
    \item[$(ii)$]\textbf{Korn's inequality.}
    \[
     \|\nabla v\|_{L^2(\Omega_{\e})}
    	\leq C\|e(v)\|_{L^2(\Omega_{\e})} \quad (v\in H_{\partial\Omega}(\Omega_\e)^3).
    \]
\end{itemize}
\end{lemma}
\begin{proof}
    For $(ii)$: \cite{DL76}, and for $(iii)$ \cite[Lemma 5.3(b)]{HJ91}.
\end{proof}

\section{Assumptions and weak formulations}\label{sec:weak_form}
We start by formulating the assumptions we place on our geometry and data.
\begin{enumerate}
    \item[\textbf{(A1)}]\label{assumptionA1} The coefficients $\alpha,\rho, c, \sigma_0\in\R$ are all positive and $L\in\R$.
    The heat conductivity $K\in\R^{3\times 3}$ is symmetric and positive definite, and the stiffness tensor $\mathcal{C}\in\R^{3\times 3\times 3\times 3}$ satisfies the symmetry conditions
    \[
    \mathcal{C}_{ijkl}=\mathcal{C}_{jikl}=\mathcal{C}_{ijlk}=\mathcal{C}_{klij}\ \ (i,j,k,l=1,2,3)
    \]
    and the positivity condition $\mathcal{C}M:M\geq c|M|^2$ for some $c>0$ and for all symmetric $M\in\R^{3\times3}$.
    \item[\textbf{(A2)}]\label{assumptionA2} The volume source densities $f_\e\in L^2(S;C^{0,1}(\overline{\Omega}))^3$ and $g_\e\in L^\infty(S;C^{0,1}(\overline{\Omega}))$ satisfy
    \[
    \sup_{\e>0}\left(\|f_\e\|_{L^2(S;C^{0,1}(\overline{\Omega}))^3}+\|g_\e\|_{L^\infty(S;C^{0,1}(\overline{\Omega}))}\right)<\infty.
    \]
    In particular, their Lipschitz constants, $L_f$ and $L_g$, are uniform in $\e$.
    \item[\textbf{(A3)}]\label{assumptionA3} The initial geometry $\Gamma_{\e}$ is periodic and $C^3$-regular and, as a consequence, its curvature function $\kappa_{\e}\colon\Gamma_\e\to\R$ is $C^1$.
    \item[\textbf{(A4)}] The initial condition $\theta_{0\e}\in L^2(\Omega)$ satisfies
    \[
    \sup_{\e>0}\|\theta_{0\e}\|_{L^2(\Omega)}<\infty.
    \]
    \item[\textbf{(A5)}] There are limit functions $f\in L^2(S\times\Omega)^3$ and $g\in L^2(S\times\Omega)$ such that $f_\e\to f$ strongly in $L^2(S\times\Omega)^3$ and $g_\e\to g$ strongly in $L^2(S\times\Omega)$ for $\e\to0$.
    \item[\textbf{(A6)}] There is a limit function $\theta_0\in L^2(\Omega)$ such that $\theta_{0\e}\to\theta_0$ strongly in $L^2(\Omega)$ for $\e\to0$.
\end{enumerate}
We point out that Assumptions (\textbf{A1})--(\textbf{A4}) enable us to solve the $\e$-moving boundary problem and Assumptions (\textbf{A5}) and (\textbf{A6}) are needed purely for the limit process $\e\to0$.

\begin{remark}[Regularity assumptions]
There is room to weaken some of the above assumptions.
For example:
\begin{itemize}
    \item The volume source densities only need to be Lipschitz continuous in the tubular neighborhood $\mathcal{U}_\e\subset\Omega$ where the interface is allowed to move, see also \cite{Wiedemann}.
    \item It is also possible to consider non-linear right hand sides.
    In case of Lipschitz continuity, e.g., the results presented in \cite{Wiedemann} could easily be adapted to the problem presented here. 
    As this introduces some additional technical difficulties, we went with the linear case here.
    \item For the initial geometry $\Gamma_\e$, $C^2$-regularity is enough to make sense of the Hanzawa transformation.
    However, the Hanzawa transformation relies on the normal vector to parametrize the transformation which incurs loss of regularity.
    To stay in $C^2$, this makes it necessary to either introduce a smooth reference interface \cite{Pruss2016} or a suitable \textit{phantom geometry} \cite{Guidotti17} and thereby introduces additional local variations in the height functions making the analysis (with respect of the scale parameter) of the problem much more involved.
\end{itemize}

\end{remark}

Our concept of weak formulation corresponding to System \ref{p:full_problem_moving} is given as: 

\begin{definition}[Weak solution in moving domains]\label{def:weak_sol_moving}
We say that 
\[
(v_\e,u_\e,\theta_\e)\in L^\infty(S)^{|\mathcal{I_\e}|}\times L^2(Q_\e)^3\times L^2(Q_\e)
\]
is a weak solution with moving domains, when 

\begin{itemize}
    \item[$(i)$] $(u_\e(t),\theta_{\e}(t))\in H^1(\Omega_\e(t))^3\times H^1(\Omega_\e(t))$ for almost all $t\in S$, $\theta_\e(0)=\theta_{0,\e}$, and $\partial_t\theta_{\e}\in L^2(Q_\e)$.
    \item[$(ii)$] The domain changes are characterized via the height function $h_\e(t,k)=\int_0^tv_\e(\tau,k)\di{\tau}$ in the sense of \cref{cor:hanzawa_stationary_eps}.
    \item[$(iii)$] The interface velocity is driven by the the local temperature averages:
    \begin{equation*}
    v_{\e}(t,k)=\frac{1}{|\Gamma_{\e,k}(t)|}\int_{\Gamma_{\e,k}(t)}\theta_\e(t,x)\di{\sigma}.
    \end{equation*}
    \item[$(iv)$] The weak form for $(u_\e,\theta_\e)$ given by
    \begin{align*}
    	\int_{\Omega_\e(t)}\mathcal{C}e(u_{\e}):e(v)\di{x}
    	-\int_{\Omega_\e(t)}\alpha\theta_{\e}\dive v \di{x}
        +\e^2\int_{\Gamma_\e(t)}\kappa_\e\sigma_0n_\e\cdot v\di{\sigma}
    	&=\int_{\Omega_\e(t)}f_{\e}\cdot v\di{x},\\
    	\int_{\Omega_\e(t)}\rho c\partial_t\theta_{\e}\varphi\di{x}
    	+\int_{\Omega_\e(t)}K\nabla\theta_{\e}\cdot\nabla \varphi\di{x}
    	+\e\int_{\Gamma_\e(t)}Lv_{\e}\varphi\di{\sigma}
    	&=\int_{\Omega_\e(t)}g_{\e}\varphi\di{x}.
    \end{align*}
    holds for all $(v,\varphi)\in H^1_{\partial\Omega}(\Omega_\e(t))^3\times H^1(\Omega_\e(t))$ and almost all $t\in S$.
\end{itemize}
\end{definition}
Due to the moving domains, it is difficult to work with \cref{def:weak_sol_moving}.
Although we are able to describe the changing domains via the height function by using the Hanzawa transformation via \cref{cor:hanzawa_stationary_eps}, we still need to refer back to some fixed reference geometry to conduct the analysis of that problem.

To that end, assume that a $C^1$-function $s_\e\colon\overline{S\times\Omega}\to\overline{\Omega}$ (such that $s_\e(t,\cdot)$ is a $C^1$-diffeomorphism for all $t\in S$) describes the moving domains and set $F_\e=Ds_\e$ and $J_\e=\det F_\e$.
We introduce reference (i.e., with respect to the reference geometry) parameter functions
\begin{subequations}\label{s:ref_coefficients}
\begin{alignat}{2}
	\mathcal{A}_\e&\colon S\times\Omega\to\R^{3\times3\times3\times3},\quad&\mathcal{A}_\e B&=\frac{1}{2}\left(F_\e^{-T} B+\left(F_{\e}^{-T} B\right)^T\right),\label{s:ref_coefficients:1}\\
	\mathcal{C}_{r,\e}&\colon S\times\Omega_\e\to\R^{3\times3\times3\times3},\quad& \mathcal{C}_{r,\e}&=J_\e \mathcal{A}_\e^T\mathcal{C}_\e\mathcal{A}_\e,\\
	\alpha_{r,\e}&\colon S\times\Omega_\e\to\R^{3\times3}, \quad&\alpha_{r,\e} &=J_\e \alpha F_{\e}^{-T} ,\\
	c_{r,\e}&\colon S\times\Omega_\e\to\R, \quad& c_{r,\e} &=J_\e\rho c_{d},\\
	K_{r,\e}&\colon S\times\Omega_\e\to\R^{3\times3}, \quad& K_{r,\e} &=J_\e F_\e^{-1} K F_\e^{-T}
\end{alignat}
the transformation related functions
\begin{alignat}{2}
w_{r,\e}&\colon S\times\Omega_\e\to\R^3, \quad& w_{r,\e} &=F_\e^{-1}\partial_ts_\e,\\
v_{r,\e}&\colon S\times\Gamma_\e\to\R^3, \quad& v_{r,\e} &=J_\e L \widetilde{v}_\e,\\
H_{r,\e}&\colon S\times\Gamma_{\e}\to\R^{3\times3}, \quad& H_{r,\e} &=J_\e\sigma_0\widetilde{\kappa}_\e F_\e^{-1} 
\end{alignat}
and the transformed source densities
\begin{alignat}{2}
	f_{r,\e}&\colon S\times\Omega_\e\to\R^3, \quad& f_{r,\e} &=J_\e \widetilde{f}_\e,\label{s:ref_coefficients:10}\\
	g_{r,\e}&\colon S\times\Omega_\e\to\R,\quad& g_{r,\e} &=J_\e \widetilde{g}_\e\label{s:ref_coefficients:11}.
 \end{alignat}
 where $\widetilde{f}$ of a function $f\colon Q_\e\to\R^n$ denotes the pullback $\widetilde{f}\colon S\times\Omega_\e\to\R^n$ defined via $f(t,x)=\widetilde{f}(t,s_\e(t,x))$.
\end{subequations}
Please note, that we have two different, but related, velocities now: The normal velocity in fixed coordinates, $v_\e$, defined on the interface $\Gamma_\e$ and the velocity of the coordinate transform, $w_\e$, defined on $\Omega_\e$.
When the motion $s_\e$ is given via the Hanzawa transform in the sense of \cref{cor:hanzawa_stationary_eps}, it automatically satisfies $J_\e w_{r,\e}\cdot n_\e=v_{r,\e}$ on $S\times\Gamma_\e$.

The corresponding weak formulation for the problem in fixed coordinates is given as: 

\begin{definition}[Weak solution in initial domains]\label{def:weak_sol_reference}
We say that 
\[
(v_{\e},u_{r,\e},\theta_{r,\e})\in L^\infty(S)^{\mathcal{I}_\e}\times L^2(S;H_{\partial\Omega}^1(\Omega_\e)^3)\times L^2(S; H^1(\Omega_\e))
\]
is a weak solution with respect to the initial domains when

\begin{itemize}
    \item[$(i)$] $\theta_{r,\e}(0)=\theta_{0,\e}$, and $\partial_t(c_{r,\e}\theta_{r,\e})\in L^2(S\times\Omega_\e)$.
    \item[$(ii)$] The domain changes are characterized via the height function $h_\e(t,k)=\int_0^tv_\e(\tau,k)\di{\tau}$ in the sense of \cref{cor:hanzawa_stationary_eps} and the functions in \eqref{s:ref_coefficients} are given in terms of the corresponding transformation $s_\e$.
    \item[$(iii)$] The velocity is driven by the the local temperature averages:
    \begin{equation*}
    v_{\e}(t,k)=\frac{1}{|\Gamma_{\e,k}|}\int_{\Gamma_{\e,k}}\theta_{r,\e}(t,x)\di{\sigma}.
    \end{equation*}
    \item[$(iv)$] The weak form for $(u_{r,\e},\theta_{r,\e})$ given by
    \begin{equation*}
    	\int_{\Omega_\e}\mathcal{C}_{r,\e}e( u_{r,\e}):e(v)\di{x}
    	-\int_{\Omega_\e}\theta_{r,\e}\alpha_{r,\e}:\nabla v \di{x}
        +\e^2\int_{\Gamma_\e}H_{r,\e}n_\e\cdot v\di{\sigma}
    	=\int_{\Omega_\e}f_{r,\e}\cdot v\di{x},
     \end{equation*}
     \begin{equation*}
    	\int_{\Omega_\e}\partial_t\left(c_{r,\e}\theta_{r,\e}\right) \varphi\di{x}
    	+\int_{\Omega_\e}\left(K_{r,\e}\nabla\theta_{r,\e}+c_{r,\e}w_{r,\e}\theta_{r,\e}\right)\cdot\nabla \varphi\di{x}
    	+\e\int_{\Gamma_\e}v_{r,\e}\varphi\di{\sigma}
    	=\int_{\Omega_\e}g_{r,\e}\varphi\di{x}
    \end{equation*}
     holds for all $(v,\varphi)\in H_{\partial\Omega}^1(\Omega_\e)^3\times H^1(\Omega_\e)$ and almost all $t\in S$.
\end{itemize}
\end{definition}
We want to point out that the two different perspectives, \cref{def:weak_sol_moving} vs \cref{def:weak_sol_reference}, are equivalent given the assumed regularity via the pullback.

In other words, $(v_\e,u_\e,\theta_\e)$ is a solution in the sense of \cref{def:weak_sol_moving} if and only if $(v_\e,\widetilde{u}_\e,\widetilde{\theta}_\e)$ is a solution in the sense of \cref{def:weak_sol_reference}.


\section{Analysis of the $\e$-problem}\label{sec:interface_tracking}
In this section, we investigate the $\e$-dependent moving boundary problem and show that there is a unique solution $(v_{\e},u_{r,\e},\theta_{r,\e})$ in the sense of \cref{def:weak_sol_reference} over some possibly small but $\e$-independent time interval $(0,t^*)\subset S$. 
We start by first establishing some important results regarding the geometry changes (\cref{ssec:tracking_of_the_interface}), before we go on to transform our problem into the initial configuration (\cref{ssec:problem_in_initial_coord}).
After showing some important $\e$-uniform estimates in \cref{ssec:estimates_bounds}, we are then able to show the existence of a unique solution utilizing the contraction mapping principle (see \cref{thm:main_existence}).

The general argument outlined in this section is as follows:
\begin{itemize}
    \item[$(i)$] \textit{Local existence for prescribed velocity:} For every $M>0$, there is a positive time horizon $T_M$ such that for every velocity function $v_\e$ bounded by $M$ leads to a unique weak solution $\theta_{r,\e}$ over $(0,T_M)$.
    \item[$(ii)$] \textit{Self mapping property:} There exists a threshold value $M^*>0$ such that for all velocity functions $v_\e$ bounded by $M^*$ the corresponding solutions $\theta_{r,\e}$ in turn lead to velocities bounded by $M^*$.
    \item[$(iii)$] \textit{Fixed-point argument:} The self-mapping has a unique fixed-point which is a solution in the sense of \cref{def:weak_sol_reference}. 
    This solution exists at least on $(0,T_{M^*})$.
\end{itemize}

\subsection{Tracking of the interface and estimates}\label{ssec:tracking_of_the_interface}
As a first step, we show that prescribed normal velocities leads to a well defined motion that tracks the growing or shrinking interfaces.
To that end, we take a function $v\colon\overline{S}\times \mathcal{I}_\e\to\R$ and look at the changing geometry.
This means that in every cell $\e(Y+k)$ for $k\in\mathcal{I}_\e$ the interface $\e(\Gamma+k)$ is growing or shrinking with velocity $\e v(t,k)$ at times $t\in\overline{S}$.

\begin{lemma}[Motion problem]\label{hanzawa_normal_vel}
    Let $v\colon\overline{S}\times \mathcal{I}_\e\to\R$ be the normal velocities describing the movement of the interfaces $\Gamma_{\e,k}$ for $k\in \mathcal{I}_\e$ and assume that they satisfy
    \[
    M:=\max_{k\in \mathcal{I}_\e}\|v(k)\|_{L^\infty(S)}<\infty.
    \]
    Then, there exist a time horizon $T_M:=\nicefrac{a^*}{10M}$ where $a^*=\min\{a_1,a_2\}$  and height parameters $h_\e\colon[0,T_M]\times\mathcal{I}_\e\to(-\nicefrac{a_1}{10},\nicefrac{a_2}{10})$ such that
    \[
    \Gamma_\e(t)=\left\{\gamma+\e h_\e\left(t,\e\left[\frac{\gamma}{\e}\right]\right)n_\e(\gamma) \ : \ \gamma\in\Gamma_{\e}\right\}\quad (t\in[0,T^*]).
    \]
    Also, the motion $s_\e\colon[0,T_M]\times\overline{\Omega}\to\overline{\Omega}$ defined via
    \begin{equation}\label{def_trafo}
    s_\e(t,x)=x+\e\begin{cases}h_\e\left(t,\e\left[\frac{P_{\e}(x)}{\e}\right]\right)n_{\e}(P_{\e}(x))\chi\left(\frac{d_{\e}(x)}{\e}\right) \quad &,\ x\in U_{\e}\\ 
    0 \quad &,\ x\notin U_{\e}\end{cases}
    \end{equation}
    satisfies $s_\e(0,\cdot)=\ID$, $s_\e(t,\overline{\Omega})=\overline{\Omega}$, and $s_\e(t,\Gamma_{\e})=\Gamma_\e(t)$.
    Moreover, $s_\e(t,\cdot)\in C^2(\overline{\Omega})$ and
    \[
    \|Ds_\e(t)\|_{L^\infty(\Omega)},\quad \|Ds_\e^{-1}(t)\|_{L^\infty(\Omega)}\leq 2.
    \]
\end{lemma}
\begin{proof}
    For any $\gamma\in\Gamma_\e$, $\e\left[\frac{\gamma}{\e}\right]$ denotes the unique $k\in\mathcal{I}_\e$ such that $\gamma\in\e(\Gamma+k)$ (see also \cref{fig:unfolding}).
    As a consequence, we can alternatively characterize the evolving surface via
    \[
    \Gamma_\e(t)=\bigcup_{i\in\mathcal{I}_\e}\left\{\gamma+\e h_\e(t,k)n_\e(\gamma) \ : \ \gamma\in\e(\Gamma+k)\right\}.
    \]
    The height functions can be explicitly calculated via
    \[
    h_\e(t,\gamma)=\int_0^tv_\e(\tau,[\gamma])\di{\tau}\qquad(\gamma\in\Gamma_{\e}).
    \]
    For \cref{cor:hanzawa_stationary_eps} to hold, we need $-a_1<10h_\e(t,k)<a_2$ which is satisfied whenever $t\leq\nicefrac{a^*}{10M}$.
    As a consequence, we can apply \cref{cor:hanzawa_stationary_eps} for every $t\in[0,T_M]$.
\end{proof}

This transformation $s_\e$ describes and tracks the changing geometry with reference to
the periodic, initial geometry $\Gamma_\e$ via the normal velocities $\e v$.
With a slight abuse of notation, we introduce the inverse (at every point in time) transformation $s_{\e}^{-1}\colon [0,T_M]\times\overline{\Omega}\to \overline{\Omega}$ as the unique function satisfying $s_{\e}^{-1}(t,s_{\e}(t,x))=x$ for all $(t,x)\in[0,T_M]\times\overline{\Omega}$.
In the following, for every $t>0$, let 
\[
V_\e(t):=L^\infty(0,t)^{\mathcal{I}_\e}.
\]
We point out that every $v_\e\in V_\e(t)$ can be uniquely identified with a piecewise constant function $\hat{v}_\e\in L^\infty((0,t)\times\Omega)$ where the constant value in any cell $\e(Y+k)\subset\Omega$ is taken to be be the corresponding $k$-value of $v_\e$.
We introduce the Jacobian matrix $F_\e(t,x)=Ds_\e(t,x)$ and its determinant $J_\e(t,x)=\det F_\e(t,x)$ and obtain the following estimates:
\begin{lemma}\label{lem:transformation_est}
Let $M>0$ and $v\in V_\e(T_M)$ such that $|v|\leq M$. 
Also, let $s_\e\colon[0,T_M]\times\overline{\Omega}\to\overline{\Omega}$ be the corresponding motion function given by \cref{hanzawa_normal_vel}.
Then,
\[
\|F_\e(t)\|_{L^\infty(\Omega)}\leq1+\frac{5M}{a^*}t,\quad \|F_\e^{-1}(t)\|_{L^\infty(\Omega)}\leq1+\frac{10M}{a^*}t,\quad \|J_\e(t)\|_{L^\infty(\Omega)}\leq 1+C\frac{M}{a^*}t,
\]
where $C>0$ does not depend on $M$, $t$, and $\e$.
The inverse motion is bounded as
\[
\|s_\e^{-1}(t)-\ID\|_{L^\infty(\Omega)}\leq 2\e Mt
\]
Moreover, there is a constant $c>0$ (independent of $M$, $t$, and $\e$) such that $c\leq J_\e(t,x)$ as well as
\[
F_\e^{-1}(t,x)F_\e^{-T}(t,x)\xi\cdot\xi\geq c|\xi|^2\quad(\xi\in\R^3).
\]
Finally, 
\[
\|\partial_tF_\e(t)\|_{L^\infty(\Omega)}+\|\partial_tF_\e^{-1}(t)\|_{L^\infty(\Omega)}+\|\partial_tJ_\e(t)\|_{L^\infty(\Omega)}\leq \frac{5M}{a^*}.
\]
\end{lemma}

\begin{proof}
The Jacobian of $s_\e$ is given via $F_\e=\ID$ on $\Omega_\e\setminus U_{\e}$ and $F_\e=\ID+\e D\psi_{\e}$ on $U_{\e}$ where $\psi_{\e}$ is given as
\[
\psi_{\e}(t,x)=h_\e(t,\e\left[\e^{-1}P_{\e}(x)\right])n_\e(P_{\e}(x))\chi\left(d_{\e}(x)\right).
\]
We note that $h_\e$ is constant for any given cell and let $k\in\mathcal{I}_\e$ and $x\in\e(Y+k)$.
Calculating the spatial derivative of $\psi_\e$ in $x$, we therefore find that
\begin{equation*}
D\psi_{\e}(t,x)=h_\e(t,k)n_\e(P_{\e}(x))\nabla\left(\chi\left(\frac{d_{\e}(x)}{\e}\right)\right)^T
+h_\e(t,k)D(n_\e(P_{\e}(x)))\chi\left(\frac{d_{\e}(x)}{\e}\right)^T.
\end{equation*}
With \cref{lem:derivatives_geometry}, we arrive at
\begin{multline}\label{derivatives_psi}
D\psi_{\e}(t,x)=\e^{-1}h_\e(t,k)\chi'\left(\frac{d_{\e}(x)}{\e}\right)\left[n_\e(P_{\e}(x))\otimes n_\e(P_{\e}(x))\right]\\
+h_\e(t,k)\chi\left(\frac{d_{\e}(x)}{\e}\right)L_{\e}(P_{\e}(x))M_{\e}(d_{\e}(x),P_{\e}(x))\left[\ID-n_\e(P_{\e}(x))\otimes n_\e(P_{\e}(x))\right].
\end{multline}
Based on the scaling properties of $\Gamma_\e$ (in particular, $|L_{\e}(\gamma)|\leq \nicefrac{1}{2a^*}$ and $|M_{\e}(r,\gamma)|\leq 2$ for all $-a_2\leq r\leq a_1$, $\gamma\in\Gamma_{\e}$; see \cref{ssec:Hanzawa} and specifically \cref{lem:derivatives_geometry})  and the properties of $\chi$, we see that
\[
|\e D\psi_{\e}(t,x)|\leq \frac{5}{a^*}\max_{k\in\mathcal{I}_\e}|h_\e(t,k)|.
\]
With $h_\e(t,k)=\int_0^tv_\e(\tau,k)\di{\tau}$ and $|v_\e|\leq M$, the first estimate
\[
\|F_\e(t)\|_{L^\infty(\Omega)}\leq1+\frac{5M}{a^*}t
\]
follows. 
With the same argument, it also follows that $|\ID-F_\e(t,x)|\leq\frac{5M}{a^*}t\leq1$.
For the inverse Jacobian, we can then employ the Neumann series to get
\[
|F_\e^{-1}(t,x)-\ID|\leq\left(1-|\ID-F_\e(t,x)|\right)^{-1}
\]
Now, as $\nicefrac{5M}{a^*}t\leq\nicefrac{1}{2}$ for all $t\in [0,T_M]$, it further follows that
\[
|F_\e^{-1}(t,x)-\ID|\leq\left(1-|\ID-F_\e(t,x)|\right)^{-1}\leq \left(1-\frac{5M}{a^*}t\right)^{-1}\leq1+\frac{10M}{a^*}t.
\]
Let $\sigma_{1,\e}(t,x)\geq\sigma_{2,\e}(t,x)\geq\sigma_{3,\e}(t,x)>0$ denote the eigenvalue functions of $F_\e$.
Then, $|F_\e(t,x)|=\sigma_{1,\e}(t,x)$,  $|F_\e^{-1}(t,x)|=\sigma_{3,\e}^{-1}$ and $|J_\e(t,x)|=\prod_{i=1}^3\sigma_{i,\e}(t,x)$.
Now, with
\[|J_\e(t,x)|\leq \sigma_{1,\e}(t,x)^3\leq\left(1+\frac{10M}{a^*}t\right)^3\leq 1+C\frac{5M}{a^*}t\]
for a sufficiently large $C>0$ (but independent of $M$, $t$, $\e$).
Similarly, we get
\begin{equation}\label{eq:estimate_determinant}
|J_\e(t,x)|\geq 1-C\frac{5M}{a^*}t\geq c>0
\end{equation}
where $c>0$ does also not depend on $M$, $t$, $\e$).
For the inverse motion, let $y\in\overline{\Omega}$ and take $x=(s_\e)^{-1}(t,y)$; this implies $y=s_\e(t,x)$ and therefore
\[
y=x+\e h_\e(t,k)n_\e(P_{\e}(x))\chi\left(d_{\e}(x)\right).
\]
Via implicit differentiation, we find that
\[
0=\partial_t(s_\e(t,s_\e^{-1}(t,y)))=\partial_ts_\e(t,s_\e^{-1}(t,y))+F_\e(t,s_\e^{-1}(t,y))\partial_ts_\e^{-1}(t,y)
\]
implying
\[
\partial_ts_\e^{-1}(t,y)=-\left[F_\e(t,s_\e^{-1}(t,y))\right]^{-1}\partial_ts_\e^1(t,s_\e^{-1}(t,y)).
\]
Here, $\partial_ts_\e(t,s_\e^{-1}(t,y))=\e v(t,k)$ where $k\in\mathcal{I}_\e$ such that $y\in\e(Y+k)$.
Integrating with respect to time, we get
\[
s_\e^{-1}(t,y)-\ID=-\e\int_0^t\left[F_\e(\tau,s_\e^{-1}(\tau,y))\right]^{-1}v(\tau,k)\di{\tau}
\]
which then yields the estimate
\[
\|s_\e^{-1}(t)-\ID\|_{L^\infty(\Omega)}\leq2\e\int_0^t\max_{k\in\mathcal{I}_\e}|v(\tau)|\di{\tau}
\]

By construction, the matrix $F^{-1}(t,x)F^{-T}(t,x)$ is symmetric and positive definite at every point $(t,x)$.
That the bound is uniform can be seen via the determinant estimates (using \cref{eq:estimate_determinant} and $\frac{5M}{a^*}t\leq\frac{1}{2}$ for all $t\in[0,T_M]$)
\[
\det(F^{-1}(t,x)F^{-T}(t,x))=\frac{1}{J_\e(t,x)^{2}}\geq\frac{1}{(1+C)^2}
\]

For the time derivatives, we first note that $v(\cdot,k)\in L^\infty(0,T_M)$ implies $h_\e(\cdot,k)\in W^{1,\infty}(0,T_M)$ for all $k\in\mathcal{I}_\e$ which in turn implies $F_\e(\cdot,x), J_\e(\cdot,x)\in W^{1,\infty}(0,T_M)$ for all $x\in\overline{\Omega}$.
With this in mind, we note that $\partial_tF_\e=0$ on $\Omega_\e\setminus U_{\e}$ and $\partial_tF_\e=\e\partial_tD\psi_{\e}$ on $U_{\e}$.
Looking at \cref{derivatives_psi}, we see that the only time dependency is due to $h_\e$.
Therefore,
\begin{multline*}
\partial_tF_\e(t,x)=\partial_th_\e(t,k)\chi'\left(\frac{d_{\e}(x)}{\e}\right)\left[n_\e(P_{\e}(x))\otimes n_\e(P_{\e}(x))\right]\\
+\partial_th_\e(t,k)\chi\left(\frac{d_{\e}(x)}{\e}\right)L_{\e}(P_{\e}(x))M_{\e}(d_{\e}(x),P_{\e}(x))\left[\ID-n_\e(P_{\e}(x))\otimes n_\e(P_{\e}(x))\right].
\end{multline*}
Now, all parts can be estimated exactly as before with the only difference being the height function.
With $h_\e(t,k)=\int_0^tv_\e(\tau,k)\di{\tau}$ and $|v_\e|\leq M$, we have $|\partial_th_\e|\leq M$.
Then (with $|\chi'|\leq\nicefrac{4m}{a^*}$, $|L_\e|\leq \nicefrac{1}{2a^*}$, and $|M_\e|\leq 2$),
\[
|\partial_tF_\e(t,x)|\leq \frac{5M}{a^*}.
\]
\end{proof}
For the fixed-point argument it is also crucial to understand how this transformation behaves under small changes in the normal velocity:

\begin{lemma}\label{lem:difference_transform}
Let $M>0$ and $v^{(1)},v^{(2)}\in V_\e(T_M)$ such that $|v^{(i)}|\leq M$.
Also, let $s_{\e}^{(i)}$ $i=1,2$, be the corresponding transformations, and let $t\in(0,T_M)$.
Then there is a constant $C>0$ (independent of $\e$, $t$, and $M$) such that
\begin{equation*}
\|s_{\e}^{(1)}-s_{\e}^{(2)}\|_{L^\infty((0,t)\times\Omega)}+\|(s_{\e}^{(1)})^{-1}-(s_{\e}^{(2)})^{-1}\|_{L^\infty((0,t)\times\Omega)}
\leq\e C\int_0^t\max_{k\in\mathcal{I}_\e}|v^{(1)}(k,\tau)-v^{(2)}(k,\tau)|\di{\tau}.
\end{equation*}
Also, for the Jacobian matrices and determinants $F_{\e}^{(i)}, J_{\e}^{(i)}$, it holds
\begin{multline*}
\|F_{\e}^{(1)}-F_{\e}^{(2)}\|_{L^\infty((0,t)\times\Omega)}+\|(F_{\e}^{(1)})^{-1}-(F_{\e}^{(2)})^{-1}\|_{L^\infty((0,t)\times\Omega)}
+\|J_{\e}^{(1)}-J_{\e}^{(2)}\|_{L^\infty((0,t)\times\Omega)}\\
\leq C\int_0^t\max_{k\in \mathcal{I}_\e}|v^{(1)}(k,\tau)-v^{(2)}(k,\tau)|\di{\tau}.
\end{multline*}
Finally, for the time derivatives,
\[
\|\partial_tJ_{\e}^{(1)}(t)-\partial_tJ_{\e}^{(2)}(t)\|_{L^\infty(\Omega)}\leq C\max_{k\in \mathcal{I}_\e}|v^{(1)}(k,t)-v^{(2)}(k,t)|.
\]
\end{lemma}
\begin{proof}
    This can easily be ascertained directly via the definition of the transformation (see \cref{def_trafo}) by incorporating the estimates established in \cref{lem:transformation_est}.
\end{proof}

In the analysis of the moving boundary problem, we are using the cell-uniform Hanzawa transformation to transform the moving geometry into the fixed reference configuration.
The following corollary is necessary for some of the estimations for the pull-back.

\begin{corollary}\label{cor:lipschitz_pull_back}
    Let $M>0$ and let $v^i\in V_\e(T_M)$ such that $|v^i|\leq M$, $i=1,2$, with corresponding Hanzawa transformations $s_\e^i\colon[0,T_M]\times\overline{\Omega}\to\overline{\Omega}$.
    Also, let $f\colon\overline{\Omega}\to\R^n$ be Lipschitz continuous (constant $L$).
    For $\tilde{f}^i\colon[0,T_M]\times\Omega$ defined via $\tilde{f}^i(t,x):=f(s_\e^{i})^{-1}(t,x))$,
    \[
    \|\tilde{f}^1(t)-\tilde{f}^2(t)\|_{L^\infty(\Omega)}\leq\e C L\int_0^t\max_{k\in\mathcal{I}_\e}|v^{(1)}(k,\tau)-v^{(2)}(k,\tau)|\di{\tau}.
    \]
\end{corollary}
\begin{proof}
    This is an immediate consequence of \cref{lem:transformation_est} and the Lipschitz continuity.
\end{proof}

\subsection{Analysis in initial domain}\label{ssec:problem_in_initial_coord}
Using the motion function $s_\e$ for a given normal velocity function $v_\e\in V_\e(T_M)$ with bound $|v_\e|\leq M$ (\cref{hanzawa_normal_vel}), the thermo-elasticity problem can be transformed into a problem posed in the initial domain.
As a direct result of \cref{lem:transformation_est} the coefficients and data functions are all bounded uniformly in $\e$, i.e., 
\[
\|\xi_{r,\e}\|_{L^\infty((0,T_M)\times\Omega)}\leq C \ \text{for}\ \xi\in\left\{\mathcal{C}_{r,\e},\alpha_{r,\e},c_{r,\e},K_{r,\e}, \frac{v_{r,\e}}{\e}, V_{r,\e}, \e H_{r,\e}, f_{r,\e},g _{r,\e}\right\}.
\]
Moreover, there exists a constant $c>0$ (independent of $\e$ and $M$) such that, for all $(0,T_M)\times\Omega$,
\begin{equation}\label{eq:positive_definite}
  c_{r,\e}(t,x)\geq c,\quad \mathcal{C}_{r,\e}(t,x)B:B\geq c|B|^2,\quad K_{r,\e}(t,x)\xi\cdot\xi\geq c|\xi|^2
\end{equation}
for all symmetric $B\in\R^{3\times3}$ and for all $\xi\in\R^3$.
Here, the positivity of $C_{r,\e}$ is not immediately obvious:
For every $u\in H^1(\Omega_\e)^3$ with $u=0$ on the external boundary $\partial\Omega$, we have
\[
\int_{\Omega_\e}C_{r,\e}e(u):e(u)\di{x}=\int_{\Omega_\e(t)}\mathcal{C}e(\tilde{u}):e(\tilde{u})\di{x}\geq c\|e(\tilde{u})\|^2_{L^2(\Omega_\e(t))},
\]
where $\tilde{u}(t,x)=u(t,s_\e(t,x))$.
With Korn's inequality and \cref{lemma:extension_operator}, this implies after transforming back to the reference coordinates
\[
\int_{\Omega_\e}C_{r,\e}e(u):e(u)\di{x}\geq c\|\nabla u\|^2_{L^2(\Omega_\e)}\leq c\|e(u)\|_{L^2(\Omega_\e)}^2.
\]
The last inequality holds because $\|\nabla u\|^2_{L^2(\Omega_\e)}=\|e(u)\|^2_{L^2(\Omega_\e)}+\|\nicefrac{1}{2}(\nabla u-\nabla u^T)\|^2_{L^2(\Omega_\e)}$ (note that $e(u)$ and $\nicefrac{1}{2}(\nabla u-\nabla u^T)$ is an orthogonal decomposition of $\nabla u$).
As this holds for all test functions $u\in H^1(\Omega_\e)^3$, positivity of $C_{r,\e}$ follows.

\begin{lemma}\label{lem:difference_transform_coefficients}
     Let $M>0$ and let $v^j\in V_\e(T_M)$ with $|v^j|\leq M$, $j=1,2$.
     Also, let $\xi_{r,\e}^j$ be any of the pull-back coefficient functions, i.e., $\xi_{r,\e}^j\in\{\mathcal{C}_{r,\e}^j,\alpha_{r,\e}^j,c_{r,\e},K_{r,\e}\}$.
     Also, let $t\in (0,T_M)$.
     There exist a $C>0$ which is independent of $\e$, $t$, and $M$ such that
     \[
     \|\xi_{r,\e}^1-\xi_{r,\e}^2\|_{L^\infty((0,t)\times\Omega)}\leq C\int_0^t\max_{k\in\mathcal{I}_\e}|v^1(\tau)-v^2(\tau)|\di{\tau}.
     \]
     For the source terms $\xi_{r,\e}\in \{f_{r,\e},g_{r,\e}\}$, it holds
     \[
     \|\xi_{r,\e}^1-\xi_{r,\e}^2\|_{L^2((0,t)\times\Omega)}\leq C(1+\e)\int_0^t\max_{k\in \mathcal{I}_\e}|v^1(\tau)-v^2(\tau)|\di{\tau},
     \]
     where the Lipschitz constants are incorporated into $C$.
     In addition, the curvature surface stresses $H_{r,\e}^j$ satisfy
     \[
     \e\|H_{r,\e}^1-H_{r,\e}^2\|_{L^\infty((0,t)\times\Gamma_\e)}\leq C\int_0^t\max_{k\in \mathcal{I}_\e}|v^1(\tau)-v^2(\tau)|\di{\tau}.
     \]
\end{lemma}
\begin{proof}
    For the coefficient functions, this follows directly from the estimates given in \cref{lem:transformation_est,lem:difference_transform} and the definition of the coefficients via the pull-back.
    For the source terms, the additional $\e$ term is due to the Lipschitz continuity and \cref{cor:lipschitz_pull_back}.
    For the curvature term, we first estimate
    \[
     \e\|H_{r,\e}^1-H_{r,\e}^2\|_{L^\infty((0,t)\times\Omega)}\leq C\int_0^t\max_{k\in I_\e}|v^1(\tau)-v^2(\tau)|\di{\tau}+C\e\|\widetilde{\kappa}_{\e}^1-\widetilde{\kappa}_{\e}^2\|_{L^\infty((0,t)\times\Gamma_\e)}
    \]
    where we have used the curvature characterization given in \cref{cor:hanzawa_stationary_eps} and $Mt<1$.
    For the second term on the right hand side, we take the corresponding height functions $h^i(t,k)=\int_0^tv^i(\tau)\di{\tau}$ and note that the curvatures are completely characterized by the height function and the point on the interface, see \cref{cor:hanzawa_stationary_eps}.
    As a consequence,
    \begin{equation*}
    \e\|\widetilde{\kappa}_{\e}^1-\widetilde{\kappa}_{\e}^2\|_{L^\infty((0,t)\times\Gamma_\e)}
    =\e\max_{k\in\mathcal{I}_\e}\|\kappa_\e(h^1(t,k),\gamma)-\kappa_\e(h^2(t,k),\gamma)\|_{L^\infty((0,t)\times\Gamma)}
    \end{equation*}
    Concentrating on the curvatures,
    \[
    \e\kappa_\e(h^i(t,k),\gamma)=\trace\left[(\ID-h^i(t,k)L_\Gamma(\gamma))^{-1}L_\Gamma(\gamma)\right],
    \]
    we see that the right hand side does not explicitly depend on $\e$.
    Moreover, $h$ is small and the right hand side depends smoothly on the height leading to
    \[
    \e|\kappa_\e(h^1(t,k),\gamma)-\kappa_\e(h^2(t,k),\gamma)|\leq C|h^1(t,k)-h^2(t,k)|.
    \]
\end{proof}

We now look at the the problem in reference coordinates for a prescribed normal velocity, i.e., we take the weak form in \cref{def:weak_sol_reference} where $v_\e\in V_\e(T_M)$ in $(iv)$ is a prescribed function.
This problem is a linear singularly coupled, implicit PDE which been analyzed in detail and shown to admit a unique weak solution with certain $\e$-uniform estimates in more complicated situations \cite{Eden_Muntean_thermo}.
For that reason, we will not reproduce the detailed analysis here but instead only present the parts of the analysis which are necessary for the fixed-point argument.

\begin{lemma}[Existence theorem for prescribed velocity]\label{lem:reference_existence}
Let $M>0$, let $v\in V_\e(T_M)$ such that $|v|\leq M$ and let Assumptions \textbf{(A1)--(A4)} be satisfied.
Then, there exists a unique $(u_{r,\e},\theta_{r,\e})\in L^2(0,T_M;H_{\partial\Omega}^1(\Omega_\e)^3\times H^1(\Omega_\e))$ with $\partial_t\theta_{r,\e}\in L^2((0,T_M)\times\Omega_\e)$ and with $\theta_{r,\e}(0)=\theta_0$ solving the variational system $(iv)$ in \cref{def:weak_sol_reference}.
\end{lemma}
\begin{proof}
This has been shown in \cite[Theorem 3.7]{Eden_Muntean_thermo} for a two-phase thermoelasticity system with a more intricate coupling.
As we are concerned with a one-phase system without dissipation, it can easily be applied in our setting.
As our velocities are cell-uniform, i.e., constant for the interfaces in a given cell $\e(Y+k)$, they also satisfy the assumption there.
Similarly, the abstract result \cite[Chapter III, Propositions 3.2 \& 3.3]{S96} can be applied directly to get a unique temperature field $\theta_{r,\e}$.
The deformation $u_{r,\e}$ can then be established via Lax-Milgram utilizing the positivity of $C_{r,\e}$.
\end{proof}
We need specific $\e$-uniform estimates of these solutions.
In addition to the standard energy estimates, we also have to establish $\e$-uniform boundedness of the temperatures to to be able to control the growth of the domains.
\begin{lemma}[Energy estimates]\label{lem:reference_apriori}
Let $M>0$ and let $v\in V_\e(T_M)$ with $|v|\leq M$.
We have
\[
	\|\theta_{r,\e}\|_{L^\infty(0,t;L^2(\Omega_{\e}))}+\|\nabla\theta_{r,\e}\|_{L^2((0,t)\times\Omega_{\e})^3}\leq C(1+\frac{M}{a^*}t),
\]
for all $t\in (0,T_M^*)$ where $C$ is independent of $M$ and $\e$.
Moreover, 
\[
\|\partial_t(c_{r,\e}\theta_{r,\e})\|_{L^2(0,T_M;H^1(\Omega_\e)^*)} + \|\theta_{r,\e}\|_{L^\infty((0,T_M)\times\Omega_\e)}\leq C.
\]
\end{lemma}
\begin{proof}
Taking the weak form for the temperature in reference coordinates and testing with $\theta_{r,\e}$, we have
\begin{equation*}
	\int_{\Omega_\e}\partial_t\left(c_{r,\e}\theta_{r,\e}\right) \theta_{r,\e}\di{x}
	+\int_{\Omega_\e}\left(c_{r,\e}w_{r,\e}\theta_{r,\e}+K_{r,\e}\nabla\theta_{r,\e}\right)\cdot\nabla \theta_{r,\e}\di{x}\\
	+\int_{\Gamma_\e}\e v\theta_{r,\e}\di{\sigma}
	=\int_{\Omega_\e}g_{r,\e}\theta_{r,\e}\di{x}.
\end{equation*}
Using the positivity of $c_{r,\e}$ and $K_{r,\e}$ and the product rule for the time derivative nets us
\begin{multline*}
    \ddt\|\theta_{r,\e}\|_{L^2(\Omega_\e)}^2
    +\|\nabla\theta_{r,\e}\|_{L^2(\Omega_\e)}^2\\
	\leq C\int_{\Omega_\e}\left(|g_{r,\e}|+|c_{r,\e}||w_{r,\e}||\nabla\theta_{r,\e}|+|\partial_tc_{r,\e}||\theta_{r,\e}|\right)|\theta_{r,\e}|\di{x}
    +\int_{\Gamma_\e}\e |v||\theta_{r,\e}|\di{\sigma},
\end{multline*}
and, making use of the estimates for $J_\e$ and $F_\e$ (\cref{lem:transformation_est}), we are led to
\begin{multline*}
    \ddt\|\theta_{r,\e}\|_{L^2(\Omega_\e)}^2
    +\|\nabla\theta_{r,\e}\|_{L^2(\Omega_\e)}^2\\
	\leq C(1+\frac{\tilde{C}M}{a^*}t)\left(\|\theta_{r,\e}\|_{L^2(\Omega_\e)}+\|\theta_{r,\e}\|_{L^2(\Omega_\e)}\|\nabla\theta_{r,\e}\|_{L^2(\Omega_\e)}\right)\\
    +CM\|\theta_{r,\e}\|_{L^2(\Omega_\e)}^2
    +\e M(1+\frac{\tilde{C}M}{a^*}t)\|\theta_{r,\e}\|_{L^1(\Gamma_\e)}.
\end{multline*}
From here (after subsuming the gradient term on the left hand side via Young's inequality), Gronwall's inequality and the trace inequality lead to the desired energy estimates.
Using these estimates and testing with any $\phi\in H^1(\Omega_\e)$ with $\|\phi\|_{H^1(\Omega_\e)}\leq1$ yields
\[
\|\partial_t(c_{r,\e}\theta_{r,\e})\|_{L^2(0,T_M;H^1(\Omega_\e)^*)}\leq C.
\]

To show that $\theta_{r,\e}$ is bounded, we introduce for $k>0$
\[
\xi_{k,\e}=\begin{cases}
        \theta_{r,\e}-k &\text{for}\ \theta_{r,\e}\geq k\\
        0&\text{for}\ |\theta_{r,\e}|\leq k\\
        \theta_{r,\e}+k &\text{for}\ \theta_{r,\e}\leq -k
        \end{cases},\quad     A_\e(t,k)=\{(t,x)\in (0,\tau)\times\Omega_\e\ : \ |\theta_{r,\e}(\tau,x)|\geq k\}
\]
and note that $\xi_{k,\e}$ is a valid test function satisfying $\partial_t \xi_{k,\e}=\chi_{A_\e(t,k)}\partial_t \theta_{r,\e}$, $\nabla \xi_{k,\e}=\chi_{A_\e(t,k)}\nabla \theta_{r,\e}$ as well as $\xi_{k,\e}\theta_{r,\e}\geq|\xi_{k,\e}|^2$.
Moreover, $\theta_{\e,r}\in L^\infty((0,T_M)\times\Omega_\e)$ if and only if there is $k>0$ such that $\xi_{k,\e}=0$ almost everywhere.

Utilizing $\xi_{k,\e}$ as a test function, we are able to estimate
\begin{multline*}
    \frac{1}{2}\|\xi_{k,\e}(t)\|_{L^2(\Omega)}^2
	+c\|\nabla\xi_{k,\e}\|_{L^2((0,t)\times\Omega)}^2\\
	\leq\frac{1}{2}\|\partial_tc_{r,\e}\|_{L^\infty((0,t)\times\Omega)}\left(\|\xi_{k,\e}\|_{L^2((0,t)\times\Omega)}^2+\|k\xi_{k,\e}\|_{L^1((0,t)\times\Omega)}\right)\\
    +\|g_{r,\e}\|_{L^2(A_\e(t,k))}\|\xi_{k,\e}\|_{L^2((0,t)\times\Omega)}\\
    -\int_0^t\int_{\Omega_\e}c_{r,\e}w_{r,\e}\theta_{r,\e}\cdot\nabla \xi_{k,\e}\di{x}\di{\tau}
	-\int_0^t\int_{\Gamma_\e}\e v\xi_{k,\e}\di{\sigma}\di{\tau}.
\end{multline*}
For the convection term, we obtain
\begin{equation*}
    -\int_0^t\int_{\Omega_\e}c_{r,\e}w_{r,\e}\theta_{r,\e}\cdot\nabla \xi_{k,\e}\di{x}\di{\tau}
    \leq 
    \e C\left(\|\xi_{k,\e}\|_{L^2((0,t)\times\Omega)}^2+\|\nabla\xi_{k,\e}\|_{L^2((0,t)\times\Omega)}^2+k^2|A_\e(t,k)|
    \right),
\end{equation*}
and for the interface term (using the trace theorem and Hölder's and Young's inequalities)
\begin{align*}
-\int_0^t\int_{\Gamma_\e}\e V_{r,\e}\xi_{k,\e}\di{\sigma}\di{\tau}
&\leq C\left(\|\xi_{k,\e}\|_{L^1((0,t)\times\Omega_\e)}+\e\|\nabla\xi_{k,\e}\|_{L^1((0,t)\times\Omega_\e)}\right)\\
&\leq C\left(\|\xi_{k,\e}\|_{L^2((0,t)\times\Omega_\e)}^2+\e^2\|\nabla\xi_{k,\e}\|_{L^2((0,t)\times\Omega_\e)}^2+|A_\e(t,k)|\right)
\end{align*}
Putting these estimate together, we have (using $\e\ll1$ and $\|g_{r,\e}\|_{\infty}\leq k$)
\[
    \|\xi_{k,\e}(t)\|_{L^2(\Omega)}^2
	+\|\nabla\xi_{k,\e}\|_{L^2((0,t)\times\Omega)}^2
	\leq C\left(\|\xi_{k,\e}\|_{L^2((0,t)\times\Omega_\e)}^2+k^2|A_\e(t,k)|\right).
\]
This implies the existence of $C>0$ (independent of $\e$ and $M$) such that $|\theta_{r,\e}|\leq C$ almost everywhere in $(0,T_M)\times\Omega_\e$. 
Here, we refer to the classical argument in \cite[Theorem 6.1]{Ladyzenskaja} and for the non-dependence on $\e$ to \cite[Lemma 6]{Wiedemann} and \cite[Theorem 3.4]{Bhattacharya22}.

%
\end{proof}

Please note that using the $C^2$ diffeomorphism $s_\e$ given by \cref{hanzawa_normal_vel}, the push-forwards of the solution $(u_{r,\e},\theta_{r,\e})$, $u_\e(t,x)=u_{r,\e}(t,s_\e^{-1}(t,x))$ and $\theta_\e(t,x)=\theta_{r,\e}(t,s_\e^{-1}(t,x))$, solve the moving boundary problem in the sense of \cref{def:weak_sol_reference} with condition $(iii)$ replaced by the prescribed velocity $v$.

\begin{lemma}[Lipschitz estimates]\label{lem:bounds_differences_temp}
   Let $j=1,2$, $M>0$, and $v^j\in V_\e(T_M^*)$ with $|v^j|\leq M$.
   Also, let $\theta_{r,\e}^j$ be the corresponding temperature solution given via \cref{lem:reference_existence}.
   Then,
   \[
    \|\bar{\theta}_{r,\e}\|_{L^\infty(0,t;L^2(\Omega_\e))}
    +\|\bar{\theta}_{r,\e}\|_{L^2(0,t;H^1(\Omega_\e))}
    \leq C(\e)\int_0^t\max_{k\in\mathcal{I}_\e}|\bar{v}(\tau,k)|\di{\tau}
   \]
   where $C(\e)$ may depend on $\e$ and where we use the bar-notation to denote differences with respect to $j$, e.g., $\bar{\theta}_{r,\e}=\theta_{r,\e}^1-\theta_{r,\e}^2$.
\end{lemma}
\begin{proof}
We test the weak form of the temperature equation with $\p=\bar{\theta}_{r,\e}$, take the difference of the two equations, and use the available estimates.
As we do not have to track the $\e$-dependicies here, this is straightforward (we refer to \cite[Lemma 6]{Wiedemann23} for details).
\end{proof}

\subsection{Fixed-point argument}\label{ssec:estimates_bounds}
In this section, we introduce a fixed-point problem formulation and combine the results from the previous \cref{ssec:tracking_of_the_interface,ssec:problem_in_initial_coord} to establish the existence of a unique solution to the System \eqref{p:full_problem_moving}.
In a first step, we utilize the results from the preceding sections to introduce the fixed-point formulation.
We emphasize that the following two critical results from \cref{ssec:problem_in_initial_coord}, \cref{lem:reference_existence,lem:reference_apriori} are valid for any $M>0$:
\begin{itemize}
    \item[$(i)$] Starting with a velocity $v\in V_\e(T_M)$ with $|v|\leq M$, we get a temperature distribution $\theta_{r,\e}\in L^2(0,T_M;H^1(\Omega_\e))$.
    \item[$(ii)$] Moreover, $\|\theta_{r,\e}\|_{L^\infty((0,T_M)\times\Omega_\e)}\leq C$ where $C>0$ does not depend on $\e$ or $M$.
\end{itemize}
This, crucially, allows us to play with the value of $M$ and show that, as long as $M$ is sufficiently large, the new velocity given by (\cref{p:full_problem_moving:4})
\begin{equation}\label{eq:new_temperature}
    \tilde{v}_\e(t,k)=\frac{1}{|\Gamma_{\e,k}|}\int_{\Gamma_{\e,k}}\theta_{r,\e}(t,x)\di{\sigma}.
\end{equation}%
is itself bounded by $M$ over $(0,T_M)$.
However, please note that larger values of $M$ imply a smaller time interval $(0,T_M)$, so choosing larger values also comes with a cost.
Combining the two statements $(i)$ and $(ii)$, we get the $M$ independent number
\begin{equation}\label{eq:definition_M}
M^*:=\sup\{\|\theta_{r,\e}\|_{L^\infty((0,T_M)\times\Omega_\e)}\ : v\in V_\e(T_M),\ |v|\leq M\ \}.
\end{equation}
Since $M^*$ is independent of $M$ it, in particular, has the same value for $M=M^*$.
Now, we set $T^*:=T_{M^*}$ and introduce, for $t\in(0,T^*)$, 
\[
\mathcal{W}_\e(t):=\left\{\vartheta\in L^2(0,t;H^1(\Omega_\e))\ : \ \|\vartheta\|_{L^\infty((0,t)\times\Omega)}\leq M^*\right\}
\]
equipped with the $L^2(0,t;H^1(\Omega_\e))$-norm.
As a closed subset of the Banach space $L^2(0,T^*;H^1(\Omega_\e))$, $\mathcal{W}_\e(t)$ is a complete metric space.

\begin{lemma}\label{lem:fixed_point_vel}
    Let $\vartheta\in \mathcal{W}_\e(T^*)$ and take the corresponding velocity according to \cref{eq:new_temperature}, i.e., 
    \[
    \tilde{v}_\e(t,k)=\frac{1}{|\Gamma_{\e,k}|}\int_{\Gamma_{\e,k}}\vartheta(t,x)\di{\sigma} \quad (t\in(0,T^*),\ k\in\mathcal{I}_\e).
    \]
    Then,  $|\tilde{v}_\e|\leq M$.
\end{lemma}

\begin{proof}
As $\vartheta(t)\in H^1(\Omega_\e)$ for almost all $t\in (0,T^*)$ as well as $\|\vartheta\|_{L^\infty((0,T^*)\times\Omega)}\leq M^*$, this follows via \cref{trace_operators}.(iii)  and \cref{eq:definition_M}.
\end{proof}
Then, starting with any $\vartheta\in\mathcal{W}_\e(t)$ for $t\in (0,T^*)$ , constructing a velocity via \cref{lem:reference_existence}, and calculating the new temperature via \cref{ssec:problem_in_initial_coord}, we land again in $\mathcal{W}_\e(t)$.
This holds for all $t\in T^*$.
More formally, \cref{eq:new_temperature} induces an operator 
\[
\mathcal{L}_\e^1\colon W_\e(t)\to\{v\in V_\e(M^*)\ : \ |v|\leq M^*\},\quad  \mathcal{L}_\e^1\vartheta(t,k)=\frac{1}{|\Gamma_{\e,k}|}\int_{\Gamma_{\e,k}}\vartheta(t,x)\di{\sigma}
\]
and, in turn, \cref{ssec:problem_in_initial_coord} induces an operator 
$\mathcal{L}_\e^2\colon\{v\in V_\e(M^*)\ : \ |v|\leq M^*\}\to  W_\e(t)$.
As a consequence, $\mathcal{L}_\e:=\mathcal{L}_\e^2\circ\mathcal{L}_\e^1\colon W_\e(t)\to W_\e(t)$ and, by construction, any function $\vartheta$ satisfying $\mathcal{L}_\e(\vartheta)=\vartheta$ solves the moving boundary problem over $(0,t)$.

In a first step, we are going to show that there is a unique weak solution for a possibly very small and $\e$-dependent time interval $(0,t_\e)$ to then argue that there must be a minimal time horizon $t^*$, independent of $\e$, such that $t^*<t_\e$.
As a result, we have a short-time solution on a time-independent interval $(0,t^*)$.
\begin{theorem}[Existence for Moving Boundary Problem]\label{thm:main_existence}
There is $t^*>0$ such that there is a unique solution $(\theta_{r,\e},v_\e)\in \mathcal{W}_\e(t^*)\times V_\e(t^*)$ to the moving boundary problem. 
\end{theorem}
\begin{proof}
We prove this statement in two steps.

 $(i)$ \textit{Existence for short $\e$-dependent time interval.}
    We need to show that $\mathcal{L}_\e$ is a contraction for a potentially small time interval $(0,T_\e)\subset(0,T^*)$, i.e., there exists a $q_\e\in(0,1)$
    \[
    \|\mathcal{L}_\e(\vartheta_1)-\mathcal{L}_\e(\vartheta_2)\|_{L^2((0,T_\e);H^1(\Omega))}\leq q_\e\|\vartheta_1-\vartheta_2\|_{L^2((0,T_\e);H^1(\Omega))}
    \]
    for all $\vartheta_1,\vartheta_2\in\mathcal{W}_\e$.

    Now, let $\vartheta_i\in \mathcal{W}_\e(t)$ and set $v_{\e,i}=\mathcal{L}_{\e}^2(\vartheta_i)$ for $i=1,2$.
    With \cref{lem:bounds_differences_temp}, we already know that
    \[
    \|\mathcal{L}_\e(\vartheta_1)-\mathcal{L}_\e(\vartheta_2)\|_{L^2(0,t;H^1(\Omega))}
    \leq C(\e)\int_0^t\max_{k\in\mathcal{I}_\e}\big|v_{\e,1}(\tau,k)-v_{\e,2}(\tau,k)\big|\di{\tau}
    \]
    where $C(\e)>0$ does not depend on $t$ or $\vartheta_i$.
    The velocities are given by
    \[
    v_{\e,i}(\tau,k)=\frac{1}{|\Gamma_{\e,k}|}\int_{\Gamma_{\e,k}}\vartheta_i(\tau,\gamma)\di{\sigma}
    \]
    and, for any $k\in\mathcal{I_\e}$ and almost every $\tau\in (0,T^*)$, it holds
    \[
    |v_{\e,1}(\tau,k)-v_{\e,2}(\tau,k)|\leq\frac{1}{|\Gamma_{\e,k}|}\int_{\Gamma_{\e,k}}|\vartheta_1(\tau,\gamma)-\vartheta_2(\tau,\gamma)|\di{\sigma}.
    \]
    With $|\Gamma_{\e,k}|=\e^2|\Gamma|$ and the trace theorem for an individual cell (\cref{trace_operators}.(ii)), this implies
    \[
    |v_{\e,1}(\tau,k)-v_{\e,2}(\tau,k)|\leq\frac{C}{\e^3|\Gamma|}\left(\|\vartheta_1(\tau)-\vartheta_2(\tau)\|_{L^1(\Gamma_{\e,k})}+\e\|\nabla\vartheta_1(\tau)-\nabla\vartheta_2(\tau)\|_{L^1(\Gamma_{\e,k})}\right)
    \]
    and therefore also
    \begin{align*}
    \max_{k\in\mathcal{I}_\e}\big|v_{\e,1}(\tau,k)-v_{\e,2}(\tau,k)\big|
    &\leq\frac{C}{\e^3|\Gamma|}\sum_{k\in\mathcal{I}_\e}\left(\|\vartheta_1(\tau)-\vartheta_2(\tau)\|_{L^1(\Gamma_{\e,k})}+\e\|\nabla\vartheta_1(\tau)-\nabla\vartheta_2(\tau)\|_{L^1(\Gamma_{\e,k})}\right)\\
    &=\frac{C}{\e^3|\Gamma|}\left(\|\vartheta_1(\tau)-\vartheta_2(\tau)\|_{L^1(\Omega_\e)}+\e\|\nabla\vartheta_1(\tau)-\nabla\vartheta_2(\tau)\|_{L^1(\Omega_\e)}\right)
    \end{align*}
    This results in
    \begin{equation*}
    \|\mathcal{L}_\e(\vartheta_1)-\mathcal{L}_\e(\vartheta_2)\|_{L^2(0,t;H^1(\Omega))}
    \leq\frac{C(\e)}{\e^3|\Gamma|}\int_0^t\left(\|\vartheta_1(\tau)-\vartheta_2(\tau)\|_{L^1(\Omega_\e)}+\e\|\nabla\vartheta_1(\tau)-\nabla\vartheta_2(\tau)\|_{L^1(\Omega_\e)}\right).
    \end{equation*}
    Utilizing Hölder's inequality in both space and time, we get 
    \[
    \|\mathcal{L}_\e(\vartheta_1)-\mathcal{L}_\e(\vartheta_2)\|_{L^2(0,t;H^1(\Omega_\e))}
    \leq\frac{C(\e)|\Omega|t}{\e^3|\Gamma|}\|\vartheta_1-\vartheta_2\|_{L^2(0,t;H^1(\Omega_\e))}.
    \]
    Consequently, $\mathcal{L}_\e$ is a contractive map on $\mathcal{W}_\e(t)$ whenever $t<\frac{\e^3|\Gamma|}{C(\e)|\Omega|}$.
    In the following, we take $T_\e=\frac{\e^3|\Gamma|}{2C(\e)|\Omega|}$ and apply the contraction mapping principle to get a unique solution $\theta_{r,\e}\in\mathcal{W}_\e(T_\e)$ satisfying $\mathcal{L}_\e(\theta_{r,\e})=\theta_{r,\e}$ over $(0,T_\e)$.

    $(ii)$ \textit{Extension of the solution to $\e$-independent interval.} At this moment, existence is only guaranteed for an $\e$-dependent time interval with $\lim_{\e\to0} T_\e=0$.\footnote{This is problematic as a homogenization limit for $\e\to0$ is meaningless in this case.}
    We want to extend the solution as far as possible and assume that $T_\e^*$ is the maximal time horizon.
    To that end, let $C_\e>0$ be the smallest constant such that
    \[
    \|\mathcal{L}_\e(\vartheta_1)-\mathcal{L}_\e(\vartheta_2)\|_{L^2(0,t;H^1(\Omega_\e))}
    \leq C_\e t\|\vartheta_1-\vartheta_2\|_{L^2(0,t;H^1(\Omega_\e))}
    \]
    for all $t\in(0,T^*_\e)$ and all $\vartheta_i\in\mathcal{W}_\e(T^*_\e)$.
    Due to the analysis in Step $(i)$, such $C_\e$ exists with $C_\e\leq\frac{C(\e)|\Omega|}{\e^3|\Gamma|}$.
    If we can establish $\sup_{\e>0}C_\e<\infty$, the statement follows.
    We take $\theta_{r,\e}\in\mathcal{W}_\e(t_{\e,\delta})$ to be the unique solution where $t_{\e,\delta}=C_\e^{-1}-\delta$ and where $\delta>0$ is a small number.
    Based on our a priori estimates, we know
    \[
   \sup_{\delta>0}\|\theta_{r,\e}\|_{L^\infty((0,t_{\e,\delta})\times\Omega_\e)}\leq M^*.
    \]
    and, because there is no norm blow-up in this time interval, the solution can be extended to $\theta_{r,\e}\in\mathcal{W_\e}(C_\e^{-1})$.
    In this situation, we have for the corresponding height functions 
    \[
    |h_\e(t,k)|\leq\int_0^t|v_\e(t,k)|\di{x}\leq C_{\e}^{-1} M^*.
    \]
    If $C_{\e}^{-1}\to0$ (even for a subsequence), we can find sufficiently small $\e_n$ for all $n\in\N$ such that $|h_\e|\leq \nicefrac{a^*}{n}$.
    This, however, implies that we can further extend the solution which in turn means that $\sup_{\e>0}C_\e<\infty$.
\end{proof}

\section{Homogenization}\label{sec:homogenization}
In this section, we present the homogenization limit of the moving boundary problem.
As similar problems and their homogenization have already be considered in detail, we skip most of the details.
For details regarding the method of two-scale convergence in general, we point to \cite{Al92,L02}, and for specific scenarios with similar homogenization procedures, to \cite{eden_homogenization_2019,Gahn21,Wiedemann}.
We use $v_\e\twosc v$ to denote two-scale convergence in the sense of \cite[Definition 1]{L02} where we say that $w_\e\in L^2(\Omega)$ two-scale converges to $w\in L^2(\Omega\times Y)$ when
\[
\int_{\Omega}w_\e(t,x)\varphi\left(x,\frac{x}{\e}\right)\di{x}
=\int_{\Omega}\int_{Y}w(t,x,y)\varphi(x,y)\di{y}\di{x}
\]
for all test functions $\varphi\in C_0^\infty(\Omega;C_\#^\infty(Y)).$
Here,
\[
C_\#^\infty(Y):=\left\{\phi\in C^\infty(\R)\ : \ \phi\ \text{is $Y$-periodic.}\right\}
\]
If, in addition,
\[
\lim_{\e\to0}\|w_\e\|_{L^2(\Omega)}\to\|w\|_{L^2(\Omega\times Y)}
\]
we say that $w_\e$ strongly two-scale converges to $w$ (notation $w_\e\stwosc w$).
If $w_\e\twosc w$ and $v_\e\stwosc v$ and $v_\e\in L^\infty(\Omega)$, then $w_\e v_\e\twosc wv$.
Please note that these definitions can easily be extended to the time-dependent setting (as outlined in \cite{L02}) treating it as a parameter.
This convergence result for products is what allows us to pass to the limit in all transformation-related quantities.
We also introduce (for $A=Y,Y_0$)
\[
H_\#^1(A):=\left\{\phi\in H_{loc}^1(\R)\ : \ \phi\ \text{is $A$-periodic}\right\}.
\]
For functions $\phi=\phi(x,y)$ depending on two spatial domains, $x\in\Omega$ and $y\in Y_0$, we use subscripts $x,y$ to signify the corresponding differential operators, i.e., $\nabla_ x\phi$ and $\nabla_y\phi$ denote the gradient with respect to $x\in\Omega$ and $y\in Y_0$, respectively.
Based on the results of the preceding section, we have a unique solution
\[
(u_{r,\e},\theta_{r,\e},v_\e)\in L^2(S^*;H_{\partial\Omega}^1(\Omega_\e))\times L^2(S^*;H^1(\Omega))\times V_\e(T^*)
\]
of the moving boundary problem and want to find the limit problem for $\e\to0$.
Here, we have set $S^*=(0,T^*)$.
For the velocities, we work with its piecewise-constant (in space) representation $\hat{v}_\e\in L^\infty(S^*\times\Omega)$.
We start by looking at the limit behavior of  functions related to the transformation:
\begin{lemma}\label{lem:limit_transformation}
   Let $\hat{v}_\e \to v$ strongly in $L^2(S^*\times\Omega)$ and $\|v\|_\infty\leq M$.
   Then the transformation  $s(t,x,\cdot)\colon\overline{Y}\to\overline{Y}$ for almost all $(t,x)\in S^*\times\Omega$ via \cref{lem:hanzawa_stationary} with height functions $h(t,x)=\int_0^tv(\tau,x)\di{\tau}$ is well defined.   
   Moreover, for the corresponding Jacobian matrix $F=Ds$ and the determinant $J=\det F$, it holds
    \begin{alignat*}{2}
        F_\e&\stwosc F,\qquad F_\e^{-1}&\stwosc F^{-1},\qquad
        J_\e&\stwosc J,\qquad \partial_tJ_\e\stwosc\partial_tJ.
    \end{alignat*}
\end{lemma}
\begin{proof}
    Since $v$ is bounded by $M$ and $S^*$ is sufficiently small, \cref{lem:hanzawa_stationary} can be applied to get the Hanzawa transformation.
    In this setup, the strong two-scale convergence was shown in \cite[Lemma 3.3]{Wiedemann23} in the stationary setting.
    For the time dependent case, see \cite[Lemma 10]{Wiedemann}.
\end{proof}
In complete analogy to the $\e$-case in \cref{s:ref_coefficients}, we introduce the limit coefficients with respect to the initial domain:
\begin{alignat*}{2}
	\mathcal{A}&\colon S\times\Omega\times Y_0\to\R^{3\times3\times3\times3},\quad&\mathcal{A} B&=\frac{1}{2}\left(F^{-T} B+\left(F^{-T} B\right)^T\right),\\
	\mathcal{C}_{r}&\colon S\times\Omega\times Y_0\to\R^{3\times3\times3\times3},\quad& \mathcal{C}_{r}&=J\mathcal{A}^T\mathcal{C}\mathcal{A},\\
	\alpha_{r}&\colon S\times\Omega\times Y_0\to\R^{3\times3}, \quad&\alpha_{r} &=J \alpha F^{-T} ,\\
	c_{r}&\colon S\times\Omega\times Y_0\to\R, \quad& c_{r} &=J\rho c_{d},\\
	K_{r}&\colon S\times\Omega\times Y_0\to\R^{3\times3}, \quad& K_{r} &=JF^{-1} K F^{-T},\\
    v_{r}&\colon S\times\Omega\times\Gamma\to\R^3, \quad& v_{r} &=JL \widetilde{v},\\
    H_{r}&\colon S\times\Omega\times\Gamma\to\R^{3\times3}, \quad& H_{r} &=J\sigma_0\widetilde{\kappa} F^{-1}. 
\end{alignat*}
With this, we can immediately follow the corresponding convergences for the coefficient functions:
\begin{corollary}
     Let $\hat{v}_\e \to v$ strongly in $L^2(S^*\times\Omega)$ and $\|v\|_\infty\leq M$.
     Then,
     \begin{alignat*}{3}
         C_{r,\e}&\stwosc C_r,&\qquad\alpha_{r,\e}&\stwosc \alpha_r,&\qquad c_{r,\e}&\stwosc c_r,\\
         K_{r,\e}&\stwosc K_r,&\qquad w_{r,\e}&\stwosc0.
     \end{alignat*}
     Moreover, it holds
    \begin{align*}
    \e^2\int_{\Gamma_\e}H_{r,\e}n_\e\cdot \psi\di{\sigma}\to\int_\Omega\int_{\Gamma}H_r n\cdot \psi\di{\sigma}\di{x},\\
    \e\int_{\Gamma_\e}v_{r,\e}\varphi\di{\sigma}\to\int_\Omega\int_\Gamma v_r\varphi\di{\sigma}\di{x}
    \end{align*}
    for all for all $\psi\in C_0^\infty(\Omega)$ and all $\varphi\in C^\infty(\overline{\Omega})$.
\end{corollary}
\begin{proof}
 With  $F_\e\stwosc F$, $F_\e^{-1}\stwosc F^{-1}$, and $J_\e\stwosc J$ the limits for  $C_{r,\e}$, $\alpha_r$, $c_{r,\e}$, and $ K_{r,\e}$ follow immediately. 
 Due to $\|\partial_ts_\e\|_\infty\leq \e C$, $w_{r,\e}\stwosc0$ is also clear.
 The surface integral limits follow from the surface integral for two-scale converging functions, see \cite[Theorem 2.1]{AllaireDamlamianHornung}.
 Please note that due to the curvature scaling, we have $H_{r,e}\sim\e^{-1}$ (see \cref{cor:hanzawa_stationary_eps}) which accounts for the additional $\e$ power in front of the surface integral.
\end{proof}

Next, we focus on the convergence of the solutions:
\begin{lemma}[Convergence limits]\label{lem:convergence_limits}
    There are functions $(u,\theta)\in L^2(S^*;H_0^1(\Omega)\times H^1(\Omega))$ such that
    \[
    \theta_{r,\e}\to\theta\quad\text{in}\ L^2(S^*\times\Omega),\quad u_{r,\e}(t)\to u(t)\quad\text{in}\ L^2(\Omega)\quad\text{for almost all $t\in S^*$}.
    \]
    In addition, there are $(u_1,\theta_1)\in L^2(S^*\times\Omega;H^1_\#(Y_0)^3\times H^1_\#(Y_0))$ such that
    \[
    \nabla u_{r,\e}\twosc\chi_{Y_0}\nabla u+\nabla_y u_1,\quad \nabla \theta_{r,\e}\twosc\chi_{Y_0}\nabla \theta+\nabla_y \theta_1.
    \]
    Also,
    \[
    \hat{v}_\e\to \theta\quad \text{in}\ L^p(S^*\times\Omega)
    \]
    for every $1\leq p<\infty$.
    Finally,
    \[
    \partial_t(c_{r,\e}\theta_{r,\e})\rightharpoonup\partial_t\left(\int_{Y_0}c_r\di{y}\,\theta\right)\quad\text{in}\ L^2(S^*\times\Omega)
    \]
    These convergences possibly only hold for subsequences.
\end{lemma}
\begin{proof}
    The weak convergence for deformation and temperature are a result of the a priori estimates (\cref{lem:reference_apriori}) and the $y$-independency of the limits a consequence of their gradient estimates.
    The pointwise in time strong convergence of $u_{r,\e}$ follows via the compact embedding $H_0^1(\Omega)\subset L^2(\Omega)$.
    The strong convergence for the temperature can be established via an compactness argument similar to \cite[Theorem 1]{Simon1986}: with the energy estimates in \cref{lem:reference_apriori} and the Lipschitz estimates in \cref{lem:difference_transform}, we can apply \cite[Proposition 6.1, Remark 6.2.(ii), Proposition 6.3.]{gahn_rigorous_2024} to get strong convergence by estimating $|J_\e(t+\Delta,x)-J_\e(t,x)|$ against $\Delta $.
    The two-scale convergence of the gradients are a direct consequence of the a priori estimates (\cite[Theorem 20]{L02}).

    For the velocity $V_\e=\e v_\e$ with $v_\e\in V_\e(T^*)$, the corresponding a priori estimate $\|\hat{v}_\e\|_{L^\infty(S^*\times\Omega)}\leq M$ only secures weak convergence in $L^p(S^*\times\Omega)$.
    With the velocity
    \[
    v_\e(t,k)=\frac{1}{\e^2|\Gamma|}\int_{\Gamma_{\e,k}}\theta_{r,\e}(t,\sigma)\di{\sigma}
    \]
    and its piecewise constant counter part $\hat{v}_\e$ (with the value of any $k\in\mathcal{I}_\e$ taken on the corresponding $\e$-cell $\e(Y+k)$), we calculate
    \[
    \int_{\Omega}\hat{v}_\e(t,x)\di{x}=\frac{1}{\e|\Gamma|}\sum_{k\in\mathcal{I_\e}}\int_{\Gamma_{\e,k}}\theta_{r,\e}(t,\sigma)\di{\sigma}=\frac{1}{\e|\Gamma|}\int_{\Gamma_{\e}}\theta_{r,\e}(t,\sigma)\di{\sigma}.
    \]
    Here, we gain one power of $\e$ as the individual cells with constant values are of size $\e$, i.e., $\e=|\e(Y+k)|$.
    With the use the periodic unfolding operator $\mathcal{T}_\e\colon L^p(S\times\Gamma_\e)\to L^p(S\times\Omega\times\Gamma)$ defined via $(\mathcal{T}_\e \varphi)(t,x,y)=\varphi(t,\e[\nicefrac{x}{\e}]+\e y)$
    and the fact $\e\int_{\Gamma_\e}\varphi\di{\sigma}=\int_\Omega\int_\Gamma \mathcal{T}_\e\varphi\di{\sigma}\di{x}$,\footnote{For a reference, we point to \cite[Section 4]{Cioranescu12}.}
    we get
    \[
    \int_{S^*}\int_\Omega\left(\theta(t,x)-\hat{v}_\e(t,x)\right)\di{x}\di{t}=\frac{1}{|\Gamma|}\int_{S^*}\int_\Omega\int_{\Gamma}\left(\theta(t,x)-\mathcal{T}_\e\theta_{r,\e}(t,x,\sigma)\right)\di{\sigma}\di{x}\di{t}.
    \]
    Now, $\mathcal{T}_\e\theta_{r,\e}\to \theta$ strongly in $L^2(S^*\times\Omega\times Y_0)$ where $\theta$ is constant in $Y_0$.
    Moreover, $\nabla_y\mathcal{T}_\e\theta_{r,\e}\to0$ since $\|\nabla\theta_{r,\e}\|_2$ is bounded uniformly in $\e$, see also \cite[Theorem 17]{Wiedemann} for more details for the same argument.
    As a result, $\hat{v}_\e\to \theta$ strongly in $L^2(S^*\times\Omega)$.

    For the time derivative we first note that there is a weak limit $\vartheta$ such that $\partial_t(c_{r,\e}\theta_{\e,r})\rightharpoonup\vartheta$ in $\ L^2(S^*\times\Omega)$ (as a consequence of the energy estimates \cref{lem:reference_apriori}) and that $c_{r,\e}\theta_{r,\e}\stwosc c_r\theta$ (due to $c_{r,\e}\in L^\infty(S^*\times\Omega)$, $c_{r,\e}\stwosc c_r$, and $\theta_{r,\e}\rightarrow \theta$ in $L^2(S^*\times\Omega)$).
    From here, we conclude that $\vartheta=\partial_t(\int_{Y_0}c_r\di{y}\theta)$ via integration by parts.
\end{proof}

Now, let $\psi\in C_0^\infty(\Omega)^3$, $\psi_1\in C_0^\infty(\Omega;C_\#^\infty(Y_0))^3$, $\phi\in C^\infty(\overline{\Omega})$, and $\phi_1\in C^\infty(\overline{\Omega};C_\#^\infty(Y_0))$ and set $\psi_\e(x)=\psi(x)+\e\psi_1(x,\nicefrac{x}{\e})$ and $\varphi_\e(x)=\varphi(x)+\e\varphi_1(x,\nicefrac{x}{\e})$.
Using $(\psi_\e,\varphi_\e)$ as test functions in the weak form \cref{def:weak_sol_moving}.(iv), we have
\begin{equation*}
    	\int_{\Omega_\e}\mathcal{C}_{r,\e}e( u_{r,\e}):e(\psi_\e)\di{x}
    	-\int_{\Omega_\e}\theta_{r,\e}\alpha_{r,\e}:\nabla \psi_\e \di{x}
        +\e^2\int_{\Gamma_\e}H_{r,\e}n_\e\cdot \psi_\e\di{\sigma}
    	=\int_{\Omega_\e}f_{r,\e}\cdot \psi_\e\di{x},
\end{equation*}
\begin{equation*}
    	\int_{\Omega_\e}\partial_t(c_{r,\e}\theta_{r,\e}) \varphi_\e\di{x}
    	+\int_{\Omega_\e}\left(K_{r,\e}\nabla\theta_{r,\e}+c_{r,\e}w_{r,\e}\theta_{r,\e}\right)\cdot\nabla \varphi_\e\di{x}
    	+\e\int_{\Gamma_\e}v_{r,\e}\varphi_\e\di{\sigma}
    	=\int_{\Omega_\e}g_{r,\e}\varphi_\e\di{x}
\end{equation*}
Using the structure of the test functions and the convergence results given in \cref{lem:convergence_limits,lem:limit_transformation} and Assumption $(A5)$, we can pass to the limit $\e\to0$ and get:
\begin{subequations}\label{eq:homlimit_raw}
\begin{multline}
    	\int_{\Omega}\int_{Y_0}\mathcal{C}_{r}\left(e( u)+e_y(u_1)\right):\left(e(\psi)+e_y(\psi_1)\right)\di{y}\di{x}
    	-\int_{\Omega}\int_{Y_0}\theta\alpha_{r}:\left(\nabla \psi+\nabla_y\psi_1\right)\di{y} \di{x}\\
        +\int_\Omega\int_{\Gamma}H_r n\cdot \psi\di{\sigma}
    	=\int_{\Omega}\int_{Y_0}Jf\cdot \psi\di{y}\di{x},
\end{multline}
\begin{multline}
    	\int_{\Omega}\partial_t\left(\int_{Y_0}c_{r}\di{y}\theta_{r}\right) \varphi\di{x}
    	+\int_{\Omega}K_{r}\left(\nabla\theta+\nabla_y\theta_{1}\right)\cdot\left(\nabla \varphi+\nabla_y\varphi_1\right)\di{y}\di{x}\\
    	+\e\int_\Omega\int_{\Gamma}v_{r}\varphi\di{\sigma}
    	=\int_{\Omega_\e}\int_{Y_0}Jg\varphi\di{x}
\end{multline}
\end{subequations}
From here there are two main steps remaining to arrive at the final limit problem: $(i).$ characterizing and eliminating $u_1$ and $\theta_1$ via solutions of cell problems and $(ii)$ going back to a formulation in moving domains.
The decoupling step is standard and it is essentially the same procedure as for other homogenization problems (cf. \cite{EKK02,Eden_Muntean_thermo,Gahn21,GahnPop23}) and, for step $(ii)$, we refer to \cite{Wiedemann23} where this was rigorously done in detail for a very general setup.

With the limit in \eqref{eq:homlimit_raw} in mind, we introduce some additional sets and functions to express the limit transformation.
First, we take $Y_0(t,x)=s(t,x,Y_0)$ and $\Gamma(t,x)=s(t,x,\Gamma)$ to denote the transformed unit cell and the interior interface, respectively.
Note that this implies $|Y_0(t,x)|=\int_{Y_0}J(t,x,y)\di{y}$.
We also introduce the cell solutions $\eta_j\in L^2(S^*\times\Omega;H^1_\#(Y_0))$, $j=1,2,3$, as the unique zero-average solution to
\[
\int_{Y_0(t,x)}K\left(\nabla_y\eta_j+e_j\right)\nabla_y\varphi\di{y}=0 \quad(\varphi\in H^1_\#(Y_0))
\]
and the functions $\mu_{jk}$, $j,k=1,2,3$, as the unique zero-average solution to
\[
\int_{Y_0(t,x)}\mathcal{C}e_y(\mu_{jk}+d_{jk}):e_y(\varphi)\di{y}=0 \quad(\varphi\in H^1_\#(Y_0)^3).
\]
Here, $e_j$ are the euclidean unit normal vectors and $d_{jk}\colon Y\to\R^3$ is given via $d_{jk}(y)=(y_j\delta_{1k},y_j\delta_{2k},y_j\delta_{3k})$ where $\delta$ is the Kronecker delta.
We further introduce the following effective coefficient functions:
\begin{alignat*}{2}
    \phi&\colon S^*\times\Omega\to\R, &\quad \phi(t,x)&=|Y_0(t,x)|,\\
    \phi_\Gamma&\colon S^*\times\Omega\to\R, &\quad \phi_\Gamma(t,x)&=|\Gamma(t,x)|,\\
    K^*&\colon S^*\times\Omega\to\R^{3\times3},&\qquad K_{ij}^*(t,x)&=\int_{Y_0(t,x)}K\left(\nabla_y\eta_j(t,x,y)+e_j\right)\cdot e_i\di{y},\\
    \mathcal{H}^*&\colon S^*\times\Omega\to\R^3,&\qquad \mathcal{H}^*(t,x)&=\int_{\Gamma(t,x)}H_\Gamma(t,x,y)n(y)\di{\sigma},\\
    \mathcal{C}^*&\colon S^*\times\Omega\to\R^{3\times3\times3\times3},&\qquad \mathcal{C}^*_{ijkl}(t,x)&=\int_{Y_0(t,x)}\mathcal{C}e_y\left(\mu_{ij}+d_{ij}\right): e_y(\mu_{kl})\di{y}.
\end{alignat*}

\begin{theorem}[Homogenization limit]\label{thm:homogenization_limit}
    The limit functions $(u,\theta)\in L^2(S^*;H_0^1(\Omega)\times H^1(\Omega))$ with $\partial_t\theta\in L^2(S^*\times\Omega)$ are the unique weak solutions of the following problem:
    \begin{subequations}\label{system:homogenized}
    \begin{equation}\label{system:homogenized:a}
    	\int_{\Omega}(\mathcal{C}^*e( u)-\theta\alpha\phi\mathds{I}):e(\psi)\di{x}
    	=\int_{\Omega}\phi f\cdot \psi\di{x}
    	+\int_{\Omega}\mathcal{H}^*\cdot \psi\di{x},
    \end{equation}
    \begin{equation}\label{system:homogenized:b}
    	\int_{\Omega}\partial_t\left(c\phi\theta\right) \varphi\di{x}\\
    	+\int_{\Omega}K^*\nabla\theta\cdot\nabla \varphi\di{x}
    	+\int_{\Omega}L\phi_\Gamma\theta\varphi\di{x}
    	=\int_{\Omega}\phi g\varphi\di{x}
    \end{equation}
    for all $(\psi,\varphi)\in H^1(\Omega)^3\times H^1(\Omega)$.
    \end{subequations}
    The microstructure evolves with the normal velocity $v=\theta$.
\end{theorem}
\begin{proof}
    With the results from \cref{lem:convergence_limits,lem:limit_transformation}, we can pass to the limit $\e\to0$ to get \eqref{eq:homlimit_raw}.
    The specific form of the limit system follows via the typical decoupling arguments using the cell solutions and the coefficient functions.
    For more details regarding this process in a similar system, we point to \cite[Section 4]{Wiedemann} and to \cite{Eden_Muntean_thermo}.
\end{proof}

Please note that this problem is still highly nonlinear, as all coefficients depend on the microstructure evolution and therefore the history of the temperature.
There is an alternative way to formulate this problem via the height functions $h(t,x)=\int_0^tv(\tau,x)\di{\tau}$.
For any given height function $h\in L^2(S^*\times\Omega)$ all relevant microstructure quantities can be calculated.
We introduce the notation $Y_0(h)$ to denote the domain corresponding to that height via \cref{lem:hanzawa_stationary} and similarly $\phi(h)$ and so on. 
The limit problem can then be written in the following way:

\begin{corollary}
    An alternative characterization of the limit problem is given as: Find $(u,\theta,h)\in L^2(S^*;H_0^1(\Omega)\times H^1(\Omega)\times L^\infty(\Omega))$ with $\partial_t\theta\in L^2(S^*\times\Omega)$ such that 
    \begin{subequations}\label{system:homogenized_alt}
    \begin{align}
    	\int_{\Omega}(\mathcal{C}^*(h)e( u)-\theta\alpha\phi(h)\mathds{I}):e(\psi)\di{x}
    	  +\int_{\Omega}\mathcal{H}^*(h)\cdot \psi\di{x}&=\int_{\Omega}\phi(h)f\cdot \psi\di{x},\\
    	\int_{\Omega}\partial_t\left(c\phi(h)\theta\right) \varphi\di{x}
    	+\int_{\Omega}K^*(h)\nabla\theta\cdot\nabla \varphi\di{x}
    	+\int_\Omega L\phi_\Gamma(h)\theta\varphi\di{x}
    	&=\int_{\Omega}\phi(h)g\varphi\di{x},\\
        h&=\int_0^t\theta\di{\tau}
    \end{align}
    for all $(\psi,\varphi)\in H^1(\Omega)^3\times H^1(\Omega)$.
    \end{subequations}
\end{corollary}

\section{Discussion}\label{sec:discussion}
We rigorously derived the homogenization limit of a phase transformation problem for a perforated media -- System \eqref{system:homogenized} is our main result. 
There are several natural ways to apply and eventually to extend our results.
For one thing, adjusting our methods for two-phase systems can be done relatively easily as long as the normal velocity is still linked to the average of the interface temperature.
The second phase would add some additional complexity due to the coupling and the homogenization result will depend on whether the heat conductivity inside the inclusions is fast or slow (corresponding to no scaling vs. scaling with $\e^2$).

Another important further question pertains to the connectedness of the inclusions: What happens in the case that the inclusions are also connected (in a similar manner as in \cite{eden2022effective})?
This naturally complicates any transformation as the non-trivial part of the transformation cannot be restricted to the individual cells.
In the current setting of cell-uniform growth this is impossible since continuity of the height function across neighboring cells is needed to construct a global transformation (see e.g, \cref{cor:hanzawa_stationary_eps}).
As a consequence, the height function would need to have the same value in every cell.
We expect this to be resolvable by smoothing between different height functions in neighboring cells.
Naturally, this leads to additional technical difficulties as the gradient of the height function now comes into play as well.

This brings us to the next important, and more difficult, extension: phases grow without information on averaged temperatures as in \cref{p:full_problem_moving:4}.
In that case, the height function is not forced to be cell-uniform and will likely vary along any individual interface $\Gamma_{\e,k}$.
Since the Hanzawa method used to transform the domains is set up to handle variable height functions, it can easily be accommodated to cope with this change.
Moreover, using the framework of maximal parabolic regularity (cf. \cite{Pruss2016}), showing the existence of a unique weak solution is feasible.
In fact, there are already promising results for a prescribed normal velocity forcing the evolution \cite{eden_homogenization_2019} or when the complete microstructure evolution is prescribed \cite{Gahn21}.
The main issue now is the availability of additional $\e$-independent regularity estimates; they are needed to ensure the construction of  an $\e$-independent interval where solutions to the microscopic problem are expected to exist; such controlled estimates are difficult to obtain for the case of perforated media.
Higher and $\e$-uniform regularity estimates (e.g. for the gradients) are needed for this to work; this is particularly difficult here due to the perforated nature of the two-phase setup; we refer to \cite{Schweizer:2000} where a specific 2D free boundary problem was treated using higher regularity estimates.
For our specific case, we do not expect that the presence of non-uniform height functions changes the structure of the homogenization limit (cf.~\cref{rem:averaging_assumption}), but in other situations this might not be the case.

Finally, we plan to take a closer look at the computability of our limit problem. 
We want to emphasize that the alternative version given by System \eqref{system:homogenized_alt} is very well suited for numerical simulations due to the availability of a precomputing strategy.
Looking first at System \eqref{system:homogenized} it is clear that although it no longer features scale heterogeneity, it is still prohibitively expensive from a computational point of view.
This is simply due to the inherent nonlinearities (moving boundary problems are posed at the microscopic level) and due to the complex scale coupling, e.g., the macroscopic temperature is driving the moving boundary at the microscopic level. 
For example, looking at the case of a Picard-type iteration resolving the nonlinearities: for every iteration, at every time step $t$, and for every node point $x$, several elliptic cell problems have to be solved.
Compare the System \eqref{system:homogenized} with System \eqref{system:homogenized_alt}, where we can set up a \textit{ precomputing strategy}: for a predefined set of height values $h_i$, calculate the corresponding functions and coefficients (e.g., $\mathcal{K}^*(h_i)$) and interpolate between these values.
As a consequence, such a procedure reduces significantly the computational cost by lowering the number of microscopic problems to be solved: 
Since the effective permeabilities/diffusivities are generally quite stable with respect to the size of inclusions (see, e.g., \cite{Ray2018}) as long as clogging and similar degeneracies are excluded, this precomputing/interpolation step should be relatively inexpensive and easy to setup.
With precomputing, we might need to solve around 10--20 cell problems overall to already get a pretty good interpolation of $h\mapsto\mathcal{K}^*(h)$ whereas without it we are talking about a number in the order of $(\Delta t\cdot H)^{-1}$ (with time step size $\Delta t$ and $H$ the macroscopic grid size).\footnote{There are alternative but related approaches like a \textit{look-up table} where already calculated values of $\mathcal{K}^*(h)$ are stored and potentially reused on the fly.}

Of course, this strategy introduces additional approximation errors and assumes some kind of stability to justify interpolation.
We are currently in the process of analyzing the needed stability property and are preparing a suitable numerical schemes for this (and similar) limit moving boundary problems.

\section*{Acknowledgement}
Michael Eden and Adrian Muntean acknowledge funding via the European Union’s Horizon Europe research and innovation program under the Marie Skłodowska-Curie fellowship project {\em{MATT}} (project nr.~101061956, \url{https://doi.org/10.3030/101061956}).

\section*{Declaration of interest statement}
On behalf of all authors, the corresponding author states that there are no conflict of interests.

\bibliographystyle{amsplain}
\bibliography{references}

\end{document}